\documentclass[12pt]{article}
\usepackage{amsfonts}
\usepackage{bm}
\usepackage{enumerate} 
\usepackage{amssymb,amsmath}
\usepackage{latexsym,mathrsfs}
\usepackage{caption}
\usepackage{subcaption}
\usepackage{accents} 
\usepackage{dsfont}
\usepackage{stackengine}
\usepackage{xr} %
\let\Horig\H

\usepackage{tikz}
\usetikzlibrary{automata,topaths}
\usetikzlibrary{shapes}
\usetikzlibrary{plotmarks}
 
\usepackage{float} 
\usepackage{color} 
\definecolor{lightblue}{rgb}{0,0.2,0.5}
\usepackage[colorlinks=true, urlcolor=lightblue, linkcolor=lightblue, citecolor=lightblue]{hyperref}

\DeclareMathAlphabet{\eufrak}{U}{}{}{} 
\SetMathAlphabet\eufrak{normal}{U}{euf}{m}{n}
\SetMathAlphabet\eufrak{bold}{U}{euf}{b}{n}

\newtheorem{assumption}{Assumption}[section]

\newcommand{\R}{\mathbb{R}}

\newcommand{\E}{\mathbb{E}}
\newcommand{\IP}{\mathbb{P}}
\newcommand{\bone}{{\bf 1}}

\newtheorem{prop}{Proposition}[section]
\newtheorem{lemma}[prop]{Lemma}
\newtheorem{definition}[prop]{Definition}
\newtheorem{corollary}[prop]{Corollary}
\newtheorem{thm}[prop]{Theorem}
\newtheorem{remark}[prop]{Remark}
\newtheorem{example}[prop]{Example}

\def\({\left(}
\def\){\right)}
\newcommand{\cov}{\mathrm{Cov}}
\def\[{\left[}
\def\]{\right]}

\def\Var{\mathrm{Var}}

\newenvironment{Proof}{\removelastskip\par\medskip
\noindent{\em Proof.} \rm}{\penalty-20\null\hfill$\square$\par\medbreak}

\allowdisplaybreaks

\numberwithin{equation}{section}

\makeatletter
\newcommand*\rel@kern[1]{\kern#1\dimexpr\macc@kerna}
\newcommand*\widebar[1]{%
  \begingroup
  \def\mathaccent##1##2{%
    \rel@kern{0.8}%
    \overline{\rel@kern{-0.8}\macc@nucleus\rel@kern{0.2}}%
    \rel@kern{-0.2}%
  }%
  \macc@depth\@ne
  \let\math@bgroup\@empty \let\math@egroup\macc@set@skewchar
  \mathsurround\z@ \frozen@everymath{\mathgroup\macc@group\relax}%
  \macc@set@skewchar\relax
  \let\mathaccentV\macc@nested@a
  \macc@nested@a\relax111{#1}%
  \endgroup
}
\makeatother

\makeatletter
\DeclareRobustCommand\widecheck[1]{{\mathpalette\@widecheck{#1}}}
\def\@widecheck#1#2{%
    \setbox\z@\hbox{\m@th$#1#2$}%
    \setbox\tw@\hbox{\m@th$#1%
       \widehat{%
          \vrule\@width\z@\@height\ht\z@
          \vrule\@height\z@\@width\wd\z@}$}%
    \dp\tw@-\ht\z@
    \@tempdima\ht\z@ \advance\@tempdima2\ht\tw@ \divide\@tempdima\thr@@
    \setbox\tw@\hbox{%
       \raise\@tempdima\hbox{\scalebox{1}[-1]{\lower\@tempdima\box
\tw@}}}%
    {\ooalign{\box\tw@ \cr \box\z@}}}
\makeatother

\usepackage{graphicx}
\usepackage{flushend,cuted}
\usepackage{bm}
\usepackage{tabularx}
\usepackage{indentfirst}
\usepackage{amssymb}
\usepackage{xparse}
\usepackage{tikz}
\usepackage{mdwlist}
\usepackage{tkz-graph}

\GraphInit[vstyle = Shade]
\usetikzlibrary[intersections,
positioning,
petri,
backgrounds,
fit,
decorations.pathmorphing,
arrows,
arrows.meta,
bending,
calc,
intersections,
through,
backgrounds,
shapes.geometric,
quotes,
matrix,
trees,
shapes.symbols,
graphs,
math,
patterns,
external,
scopes,
matrix,
lindenmayersystems,
shapes.callouts,
shapes.misc,
angles,
shapes.arrows,
shadings]
\begin{document}
\title{
\huge
 Gaussian fluctuations of generalized $U$-statistics and
 subgraph counting in the binomial random-connection model 
} 

\author{
  Qingwei Liu\footnote{School of Mathematics and Statistics, The University of Melbourne, Parkville, VIC 3010, Australia.
Email: \href{mailto:qingwei.liu@unimelb.edu.au}{qingwei.liu@unimelb.edu.au}}
  \qquad
      Nicolas Privault\footnote{Division of Mathematical Sciences, School of Physical and Mathematical Sciences, Nanyang Technological University, 21 Nanyang Link, Singapore 637371, Email: \href{mailto:nprivault@ntu.edu.sg}{nprivault@ntu.edu.sg}
}
}

\maketitle

\vspace{-0.5cm}

\begin{abstract} 
  We derive normal approximation bounds
  for generalized $U$-statistics of the form 
 \begin{equation*}
    S_{n,k}(f):=\sum_{
      1 \leq \beta (1),\dots,\beta (k) \leq n \atop
      \beta (i)\ne\beta (j), \ 1\leq i\ne j \leq k}
  f\big(X_{\beta (1)},\dots,X_{\beta (k)},Y_{\beta (1),\beta (2)},\dots,Y_{\beta (k-1),\beta (k)}\big), 
 \end{equation*}
 where $\{X_i \}_{1 \leq i \leq n}$ and $\{Y_{i,j}\}_{1\le i < j \le n}$ are
 independent sequences of i.i.d. random variables.
 Our approach relies on moment identities and cumulant bounds
 that are derived using partition diagram arguments.
 Normal approximation bounds in the Kolmogorov distance
 and moderate deviation results are then obtained
 by the cumulant method.
 Those results are applied to subgraph counting 
 in the binomial random-connection model,
 which is a generalization of the Erd{\H o}s-R\'enyi model. 
\end{abstract}

\medskip


\baselineskip0.7cm

\section{Introduction}
\noindent 
 Second-order $U$-statistics
 can be viewed as quadratic random functionals of the form
\begin{equation}
\label{fdjkl132} 
\sum_{1\leq i \ne j \leq n} Y_{i,j} X_i X_j
\end{equation} 
 which are used to model potentials and partition functions 
 in the framework of
 the Gaussian Unitary Ensemble in statistical mechanics, 
 where $(Y_{i ,j})_{1\leq i \not= j \leq n}$ is a possibly random adjacency matrix  
 made of independent entries,
 and $(X_i )_{1\leq i \leq n}$ is a sequence of independent
 identically distributed random variables. 
 Cumulant bounds of \eqref{fdjkl132} have been obtained in
 \cite{khorunzhiy} when 
 $(Y_{i ,j })_{1\leq i \not= j \leq n}$ represents the adjacency matrix
 of the Erd{\H o}s-R\'enyi random graph,
 and general approximation results in distribution have been recently
 obtained in \cite{bhattacharya2} when
 $(Y_{i ,j })_{1\leq i \not= j \leq n}$ is deterministic.
 On the other hand,
 more general pair interactions can be modeled using
 $U$-statistics of the form
\begin{equation}
\nonumber %
\sum_{1\leq i \ne j \leq n} f(X_i ,X_j ), 
\end{equation} 
while higher-order $U$-statistics can
model nonlinear interactions in non-Gaussian frameworks. 

  \medskip 

 In this paper, we derive normal approximation
 results and cumulant bounds for 
 generalized $U$-statistics of the form
\begin{equation}
  \label{def:gUstat-0}
    S_{n,k}(f):=\sum_{
      1 \leq \beta (1),\dots,\beta (k) \leq n \atop
      \beta (i)\ne\beta (j), \ \! 1\leq i\ne j \leq k}
  f\big(X_{\beta (1)},\dots,X_{\beta (k)},Y_{\beta (1),\beta (2)},\dots,Y_{\beta (1),\beta (k)},Y_{\beta (2),\beta (3)},\dots,Y_{\beta (k-1),\beta (k)}\big), 
 \end{equation}
 which have been introduced in \cite{Janson91}
 as a powerful tool for studying the normal and non-normal asymptotic distributions of subgraph counts in inhomogeneous random graphs. 
 Here, $\{X_i \}_{1\leq i \leq n}$ and $\{Y_{i ,j }\}_{1\leq i <j \leq n}$ are
 two independent sequences of i.i.d. random elements taking values
 respectively in a Borel space
$\mathcal{S}$ and a measurable space $\mathcal{M}$,
and $f:\mathcal{S}^k\times\mathcal{M}^{k(k-1)/2}\to\R$
 is a measurable function, $k\ge2$, 
with $Y_{i ,j }=Y_{j ,i }$ if $1 \leq j <i \leq n$.

\medskip
  
 In Corollary~\ref{thm4.2}, we obtain a Kolmogorov distance bound of the form 
    \begin{equation}
      \nonumber 
      \sup_{x\in\R}\left|
      \IP\left(
  \widebar{S}_{n,k}(f)
      \leq x\right)-\Phi(x)\right|\leq \frac{C(f,k)}{n^{1/(2+4k)}},
 \quad n\ge4(k-1), 
    \end{equation}
    for the normalized generalized $U$-statistics 
  $$
  \widebar{S}_{n,k}(f):=
        \frac{S_{n,k}(f)-\kappa_1(S_{n,k}(f))}{\sqrt{\kappa_2(S_{n,k}(f))}},
        $$
        where $\Phi$ is the cumulative distribution function of the standard normal distribution and $C(f,k)>0$ depends only on $f$ and on 
        $k \geq 2$. 
  Our approach relies on 
  moment identities and cumulant bounds for
    generalized $U$-statistics
    established in Theorems~\ref{jklfd12} and
    \ref{jkldd12}.       
 Berry-Esseen bounds for the normal approximation of general
functionals of binomial point processes have been obtained in 
\cite{lachieze-rey-peccati}, with application to
 $U$-statistics and set approximation 
 for random tessellations. 
However, generalized $U$-statistics of the form \eqref{def:gUstat-0}
include an additional layer of randomness 
due to the random sequence $\{Y_{i ,j }\}_{1\leq i < j \leq n}$.   
 In Corollary~\ref{corg2-0} we obtain
 a moderate deviation result for 
 the normalized $U$-statistics 
 $\widebar{S}_{n,k}(f)$. 
 
\medskip
    
  Starting with Section~\ref{s5}, we
 apply our normal approximation results
  for generalized $U$-statistics 
 to subgraph counting in the binomial random-connection model.   
Random-connection models (RCMs) are random graphs which are based on
randomly located vertices which are independently connected
with a location-dependent probability. 
 As a generalization of the Erd{\H o}s-R\'enyi model, the binomial RCM 
 has gained significant attention and has been studied under different names, for example
 as inhomogeneous random graphs,
 c.f. \cite{DevroyeFraiman14,penrose18,hladky21}, and as
 graphon-based random graphs c.f. \cite{coulson16,zhangzs,bhattacharya23}.

 \medskip

  Distributional approximations for count statistics on random-connection models
 whose vertices are generated according to
 a Poisson point process,
 have been investigated in a number of recent works, including vertices counts \cite{penrose18}, component counts \cite{LNS21}, and subgraph counts \cite{can2022,LiuPrivault}. Recently, Poisson approximation with bounds for subgraph counts
 in general random-connection models have been derived in \cite{liu-xia},
 and the cumulant method has been applied to
 subgraph counting weighted random connection models
 in \cite{heerten}.
  
\medskip
 
More formally, let
 $\mathcal{X}_n=\{X_1,\dots,X_n\}$ 
 denote a family of i.i.d. random points with a common distribution $\mu$ on $\mathcal{S}=\R^d$,
 for some $n\geq 2$. 
Given $H:\R^d\times\R^d\to[0,1]$
 a symmetric measurable connection function,
 the binomial random-connection model 
 with connection function $H$ is
 the random graph %
 on the binomial point process
$\mathcal{X}_n=\{X_1,\dots,X_n\}$
 constructed by
adding edges independently with probability $H(X_i ,X_j )$,
 to each distinct pair $(X_i ,X_j )$ of vertices,
 $1\leq i \ne j \leq n$.

 \medskip 

 When the connection function is taken
 as $H(x,y)=\bone_{\{\|x-y\|\leq r\}}$, $x,y\in\R^d$ for some $r>0$, vertices are
 connected in a deterministic way, and the binomial RCM becomes a random geometric graph, or Gilbert graph, c.f. \cite{penrosebk}. When $H(x,y)\equiv p_n$, $x,y\in\R^d$, the binomial RCM recovers the classical Erd{\H o}s-R\'enyi random graph.  In this case,  subgraph counting in the binomial RCM is a natural extension of the subgraph counting problem in the Erd{\H o}s-R\'enyi model, see \cite{JLR,feray16}.

   \medskip

   Although the asymptotic behavior of subgraph counts in the binomial RCM was studied in detail in \cite{Janson91,janson,bhattacharya23},
   convergence rates for
   the distributional approximation of subgraph counts
   have only been recently discussed, see \cite{KaurRollin21,zhangzs}. In \cite{dsc},
the Poisson approximation of standard (not generalized)
$U$-statistics has been considered in the binomial model.
More recently, general approximation results
for standard second-order $U$-statistics have been
obtained in \cite{bhattacharya2}. 
However, none of those works,
 including \cite{zhangzs}, 
   consider the case where
   the probabilities of connecting two vertices
   tends to zero as $n$ tends to infinity, 
   as is typical in the Erd{\H o}s-R\'enyi model. 

   \medskip

 In contrast with the centered subgraph counts considered in \cite{KaurRollin21},
  we allow the connection probability
   of any distinct pair $(X_i ,X_j )$ of vertices in the binomial RCM
   to be of the form $p_n H(X_i ,X_j )$
   where $p_n\in(0,1)$, $1\leq i \ne j \leq n$. 
   We consider in particular the case where $p_n$ may tend to zero
   as $n$ tends to infinity,
   and study the corresponding normal approximation of subgraph counts.
   
   \medskip

In Theorem~\ref{fjkdl243} 
 we derive upper bounds on the cumulants  
 of subgraph counts in the binomial random-connection
 model.
   Note that
  related cumulant bounds have been obtained in the Erd{\H o}s-R\'enyi model
   in \cite{khorunzhiy} for the counts of line and cycle graphs, 
   and in \cite{feray16} for general subgraph counts. 
   In comparison with the Poisson random-connection model
     considered in \cite{LiuPrivault,LiuPrivault24},
     cumulants admit no
     simplified expression using
     sums over connected non-flat partitions in the binomial
     random-connection model.
     For this reason, cumulant bounds
     have to be derived using specific
     arguments. 

   \medskip

 Then, by combining Theorem~\ref{fjkdl243}
 with variance lower bounds obtained
 in Proposition~\ref{p6.1}, in Theorem~\ref{th7.5}
 we obtain cumulant growth rates for the counts of
 strongly balanced connected graphs. 
 Our proof relies on dependency graph methods and the 
 convex analysis of planar diagrams, which were introduced in 
 \cite{luczak} to study the behaviour of the variance of
 subgraph counts in the Erd{\H o}s-R\'enyi model.

 \medskip

 Cumulant growth rates for the normalized counts of
 strongly balanced connected graphs
 are then obtained in Corollary~\ref{cjkfl}
 under Assumption~\ref{assu1} in the case where $p_n=o(1)$,
 using the variance lower bound for subgraph counts
 established in Proposition~\ref{prop7.6}. 
 Note however that Assumption~\ref{assu1} is
valid in the binomial RCM, and is
not satisfied in the Erd{\H o}s-R\'enyi model.    

\medskip
 
 Using the cumulant method and
 the Statulevi\v{c}ius condition,
 see Appendix~\ref{appx}, Kolmogorov distance rates
 of the form 
  \begin{eqnarray*}
    \sup_{x\in\R}|\IP(\widebar{N}_G\leq x)-\Phi(x)|\leq \begin{cases}
      \displaystyle
      \frac{C}{n^{1/(2+4v(G))}}~~~~~\quad~~~~{when}~~~n^{-(v(G)-1)/e(G)}\ll p_n, %
      \medskip
      \\
      \displaystyle
      \frac{C}{\big(n^{v(G)}p_n^{e(G)}\big)^{1/(2+4v(G))}}~~~~{when}~~~p_n\ll n^{-(v(G)-1)/e(G)},
    \end{cases}
  \end{eqnarray*}
  $n\geq 4(v(G)-1)$, 
  are obtained in Corollary~\ref{corg}
  for the normalized subgraph count
  $$
  \widebar{N}_G:= \frac{N_G-\kappa_1(N_G)}{\sqrt{\kappa_2(N_G)}}
  $$
  of a strongly balanced connected graph $G$ with
  $v(G)$ vertices. This extends 
  the results of \cite{PS} and \cite{rednos}
  for the approximation of subgraph counts from the
  Erd{\H o}s-R\'enyi model to the binomial RCM.
  Although the convergence rates in the Kolmogorov
    distance obtained in do not match the
    optimal rate obtained in \cite{zhangzs},
  they allow us to consider the case where $p_n = o(1)$.
  In Corollary~\ref{graphcontain}, we investigate the threshold phenomenon for the containment of subgraphs in the binomial RCM.
  
 \medskip

 In Corollary~\ref{corg2} we obtain
 moderate deviation results for 
  the normalized subgraph count
  $\widebar{N}_G$,
  see also 
  \cite[Theorem~1.1]{doring1},
  \cite[Theorem~2.3]{doring},
  \cite[Theorem~10.0.2]{feray16},
  \cite[Theorem~1.5]{alvarado}
  for moderate deviation results
  in the Erd{\H o}s-R\'enyi model. 
 In Corollary~\ref{graphcontain}, we investigate the threshold phenomenon for the containment of subgraphs in the binomial RCM.

\medskip

The paper is organized as follows. In Section~\ref{s2} we recall some notation and definitions related to set partitions and diagrams. In Section~\ref{s3} we derive moment identities for generalized $U$-statistics. Section~\ref{s4} gives cumulant bounds for generalized $U$-statistics and further obtains normal approximation results via the cumulant method. In Section~\ref{s5}, we obtain cumulant bounds for subgraph counts in the binomial random-connection model, which allows the connection probability between pairs of vertices to be significantly small. In Section~\ref{s6}, we derive lower bounds on the variance of subgraph counts in the binomial RCM. This is crucial for proving the results in Section~\ref{s7}. In Section~\ref{s7} we obtain normal approximation for subgraph counts in the binomial RCM through a refined analysis on the cumulant growth rates and a threshold phenomenon for subgraph containment. In the Appendix~\ref{appx}, we provide a brief review of the Statulevi\v{c}ius condition and its application to the cumulant method.

\section{Set partitions and diagrams}
\label{s2}
\noindent
In what follows, we let $[n]:=\{1,2,\dots,n\}$ for $n\geq 1$,
 and let $\Pi(b)$ denote the collection of set partitions of
any finite set $b$.
 Given two set partitions $\rho_1, \rho_2\in\Pi(b)$, we say that $\rho_1$ is
coarser than $\rho_2$
 (i.e. $\rho_2$ is finer than $\rho_1$), 
 and we write $\rho_2\preceq\rho_1$, 
if and only if each block of $\rho_2$ is contained in a block of $\rho_1$.
We use $\rho_1\vee\rho_2$ for the
finest partition which is coarser than both $\rho_1$ and $\rho_2$, and
denote by $\rho_1\wedge\rho_2$ the coarsest partition
which is finer than both of $\rho_1$ and $\rho_2$. We also let $\widehat{1}:=\{b\}$ 
denote the single-block coarsest partition of $b$,
 whereas $\widehat{0}$ stands for the partition made of singletons. 
 Given $k\geq 2$ and $j\geq 1$,
 we let $\pi:=\{\pi_1,\dots,\pi_j\}$ %
 denote the partition
 of $[j]\times[k]$ defined as 
 $$\pi_i:=\{(i,\ell):1\leq \ell\leq k\},
 \quad i =1,\dots , j.
 $$
 For $j,k \geq 1$ we also let $\pi_\eta := (\pi_i)_{i\in\eta}\in\Pi(\eta \times [k])$ denote the partition made of $|\eta|$ blocks of size $k$.
 A partition $\rho\in\Pi([j]\times[k])$ is said to be {\it non-flat} if $\rho\wedge\pi=\widehat{0}$,
 see Chapter~4 of \cite{peccatitaqqu} and Figure~\ref{fig:diagram0}. 
\begin{figure}[H]
    \captionsetup[subfigure]{font=footnotesize}
\centering
\subcaptionbox{Non-flat partition.}[.49\textwidth]{%
 \begin{tikzpicture}[scale=0.8] 
     \draw[black, thick] (0,0) rectangle (5,4);
 \node[anchor=east,font=\small] at (0.8,3) {1};
 \node[anchor=east,font=\small] at (0.8,2) {2};
 \node[anchor=east,font=\small] at (0.8,1) {3};
 \node[anchor=south,font=\small] at (1,0) {1};
 \node[anchor=south,font=\small] at (2,0) {2};
 \node[anchor=south,font=\small] at (3,0) {3};
 \node[anchor=south,font=\small] at (4,0) {4};
 \filldraw [black] (1,1) circle (2pt);
 \filldraw [black] (2,1) circle (2pt);
 \filldraw [black] (3,1) circle (2pt);
 \filldraw [black] (4,1) circle (2pt);
 \filldraw [black] (1,2) circle (2pt);
 \filldraw [black] (2,2) circle (2pt);
 \filldraw [black] (3,2) circle (2pt);
 \filldraw [black] (4,2) circle (2pt);
 \filldraw [black] (1,3) circle (2pt);
 \filldraw [black] (2,3) circle (2pt);
 \filldraw [black] (3,3) circle (2pt);
 \filldraw [black] (4,3) circle (2pt);
 \draw[very thick] (1,3) -- (1,2) -- (1,1);
 \draw[very thick] (3,3) -- (3,2);
 \draw[very thick] (2,2) -- (2,1);
 \draw[very thick] (3,1) -- (4,2);
  \begin{pgfonlayer}{background}
    \filldraw [line width=4mm,black!3]
      (0.2,0.2)  rectangle (4.8,3.8);
  \end{pgfonlayer}
\end{tikzpicture}
}
\subcaptionbox{Flat partition.}[.49\textwidth]{%
  \begin{tikzpicture}[scale=0.8] 
 \draw[black, thick] (0,0) rectangle (5,4);
 \node[anchor=east,font=\small] at (0.8,3) {1};
 \node[anchor=east,font=\small] at (0.8,2) {2};
 \node[anchor=east,font=\small] at (0.8,1) {3};
 \node[anchor=south,font=\small] at (1,0) {1};
 \node[anchor=south,font=\small] at (2,0) {2};
 \node[anchor=south,font=\small] at (3,0) {3};
 \node[anchor=south,font=\small] at (4,0) {4};
 \filldraw [black] (1,1) circle (2pt);
 \filldraw [black] (2,1) circle (2pt);
 \filldraw [black] (3,1) circle (2pt);
 \filldraw [black] (4,1) circle (2pt);
 \filldraw [black] (1,2) circle (2pt);
 \filldraw [black] (2,2) circle (2pt);
 \filldraw [black] (3,2) circle (2pt);
 \filldraw [black] (4,2) circle (2pt);
 \filldraw [black] (1,3) circle (2pt);
 \filldraw [black] (2,3) circle (2pt);
 \filldraw [black] (3,3) circle (2pt);
 \filldraw [black] (4,3) circle (2pt);
 \draw[very thick] (1,3) -- (1,2) -- (1,1) -- (2,1);
 \draw[very thick] (3,3) -- (3,2);
 \draw[very thick] (2,2) -- (2,1);
 \draw[very thick] (3,1) -- (4,2);
  \begin{pgfonlayer}{background}
    \filldraw [line width=4mm,black!3]
      (0.2,0.2)  rectangle (4.8,3.8);
  \end{pgfonlayer}
 \end{tikzpicture}%
}
\caption{Examples of partition diagrams with $n=3$ and $r=4$.}
 \label{fig:diagram0}
 \end{figure}
\vskip-0.3cm
\noindent
A partition $\rho\in\Pi([j]\times[k])$ is said to be {\it connected} if $\rho\vee\pi=\widehat{1}$, see Figure~\ref{fig:diagram0-1}. 
\begin{figure}[H]
\captionsetup[subfigure]{font=footnotesize}
\centering
\subcaptionbox{Connected partition.}[.49\textwidth]{
\begin{tikzpicture}[scale=0.8] 
\draw[black, thick] (0,0) rectangle (5,6);
\node[anchor=east,font=\small] at (0.8,5) {1};
\node[anchor=east,font=\small] at (0.8,4) {2};
\node[anchor=east,font=\small] at (0.8,3) {3};
\node[anchor=east,font=\small] at (0.8,2) {4};
\node[anchor=east,font=\small] at (0.8,1) {5};
\node[anchor=south,font=\small] at (1,0) {1};
\node[anchor=south,font=\small] at (2,0) {2};
\node[anchor=south,font=\small] at (3,0) {3};
\node[anchor=south,font=\small] at (4,0) {4};
\filldraw [black] (1,1) circle (2pt);
\filldraw [black] (2,1) circle (2pt);
\filldraw [black] (3,1) circle (2pt);
\filldraw [black] (4,1) circle (2pt);
\filldraw [black] (1,2) circle (2pt);
\filldraw [black] (2,2) circle (2pt);
\filldraw [black] (3,2) circle (2pt);
\filldraw [black] (4,2) circle (2pt);
\filldraw [black] (1,3) circle (2pt);
\filldraw [black] (2,3) circle (2pt);
\filldraw [black] (3,3) circle (2pt);
\filldraw [black] (4,3) circle (2pt);
\filldraw [black] (2,3) circle (2pt);
\filldraw [black] (1,4) circle (2pt);
\filldraw [black] (2,4) circle (2pt);
\filldraw [black] (3,4) circle (2pt);
\filldraw [black] (4,4) circle (2pt);
\filldraw [black] (1,5) circle (2pt);
\filldraw [black] (2,5) circle (2pt);
\filldraw [black] (3,5) circle (2pt);
\filldraw [black] (4,5) circle (2pt);
\draw[very thick] (1,5) -- (1,4); 
\draw[very thick] (2,2) .. controls (2.5,3) .. (2,4);
\draw[very thick] (2,3) .. controls (1.5,4) .. (2,5);
\draw[very thick] (3,5) -- (4,4);
\draw[very thick] (1,2) -- (1,1);
\draw[very thick] (2,1) -- (3,2) -- (4,3) -- (3,4);
  \begin{pgfonlayer}{background}
    \filldraw [line width=4mm,black!3]
      (0.2,0.2)  rectangle (4.8,5.8);
  \end{pgfonlayer}
\end{tikzpicture}}
\subcaptionbox{Non-connected partition.}[0.49\textwidth]{
\begin{tikzpicture}[scale=0.8] 
\draw[very thick,dashed,red] (0.3,3.5) -- (4.7,3.5);
\draw[black, thick] (0,0) rectangle (5,6);
\node[anchor=east,font=\small] at (0.8,5) {1};
\node[anchor=east,font=\small] at (0.8,4) {2};
\node[anchor=east,font=\small] at (0.8,3) {3};
\node[anchor=east,font=\small] at (0.8,2) {4};
\node[anchor=east,font=\small] at (0.8,1) {5};
\node[anchor=south,font=\small] at (1,0) {1};
\node[anchor=south,font=\small] at (2,0) {2};
\node[anchor=south,font=\small] at (3,0) {3};
\node[anchor=south,font=\small] at (4,0) {4};
\filldraw [black] (1,1) circle (2pt);
\filldraw [black] (2,1) circle (2pt);
\filldraw [black] (3,1) circle (2pt);
\filldraw [black] (4,1) circle (2pt);
\filldraw [black] (1,2) circle (2pt);
\filldraw [black] (2,2) circle (2pt);
\filldraw [black] (3,2) circle (2pt);
\filldraw [black] (4,2) circle (2pt);
\filldraw [black] (1,3) circle (2pt);
\filldraw [black] (2,3) circle (2pt);
\filldraw [black] (3,3) circle (2pt);
\filldraw [black] (4,3) circle (2pt);
\filldraw [black] (2,3) circle (2pt);
\filldraw [black] (1,4) circle (2pt);
\filldraw [black] (2,4) circle (2pt);
\filldraw [black] (3,4) circle (2pt);
\filldraw [black] (4,4) circle (2pt);
\filldraw [black] (1,5) circle (2pt);
\filldraw [black] (2,5) circle (2pt);
\filldraw [black] (3,5) circle (2pt);
\filldraw [black] (4,5) circle (2pt);
\draw[very thick] (1,5) -- (1,4) -- (2,4) -- (3,4);
\draw[very thick] (2,5) -- (3,5) -- (4,5) -- (4,4);
\draw[very thick] (1,2) -- (1,1);
\draw[very thick] (2,3) -- (2,2);
\draw[very thick] (2,1) -- (3,1) -- (4,1);
\draw[very thick] (3,2) -- (4,3);
  \begin{pgfonlayer}{background}
    \filldraw [line width=4mm,black!3]
      (0.2,0.2)  rectangle (4.8,5.8);
  \end{pgfonlayer}
\end{tikzpicture}}
\caption{Examples of partition diagrams with $n=5$ and $r=4$.}
\label{fig:diagram0-1}
\end{figure}

 \vspace{-.4cm}
  
\noindent 
 We also let $\Pi_{\widehat{1}}([j]\times[k])$ denote the collection of all connected partitions of $[j]\times[k]$,
 and denote by 
      $$
      {\rm CNF}(j,k) :=\big\{\rho : \rho\in\Pi_{\widehat{1}}( [j]\times [k]),\ \rho\wedge\pi=\widehat{0}\big\}
      $$
       the set of all connected and non-flat partitions of $[j]\times [k]$, for $j,k\geq 1$.   
 Let $[n]^k_{\ne}$ denote the collection of distinct $k$-fold indexes 
 \begin{equation}
   \nonumber %
  [n]^k_{\ne}:=\{\beta =(\beta (1),\dots,\beta (k))\in[n]^k : \beta (i)\ne\beta (j)~\mathrm{for}~1\leq i\ne j\leq k\},
   \qquad k \geq 1. 
\end{equation}
 The role of the partition
   $\sqcap (\alpha )$ introduced in Definition~\ref{def:parti}
   is to group each set of identical entries in a
   family of $k$-tuples into a partition block.
   Later on, it will be used to identify the
   common random variables appearing
   in repeated copies of $(X_1, \ldots , X_n)$
   for the computation of joint cumulants. 
\begin{definition}
\label{def:parti}
 Given a sequence 
 $$
 \alpha   =
\left[ 
\begin{array}{c} 
\alpha_1 
\\ 
\vdots 
\\ 
\alpha_j 
\\ 
\end{array} 
\right] 
  =
\left[ 
\begin{array}{ccc}
 \alpha_1 (1) & \cdots & \alpha_1(k)
  \\
\vdots & \ddots & \vdots 
\\
 \alpha_j (1) & \cdots & \alpha_j(k)
\end{array} 
\right] 
\in [n]^k_{\ne}\times \cdots \times [n]^k_{\ne}, 
$$ 
  we let $\sqcap (\alpha )$ denote the partition 
  of $[j]\times[k]$ such that each block of $\sqcap (\alpha )$
  is made of elements $(i,\ell )$
  that correspond to a same value of $\alpha_i (\ell )$.
\end{definition}
Next is an example of a sequence
$\alpha %
\in [n]^k_{\ne}\times \cdots \times [n]^k_{\ne}$
 and of the partition $\sqcap (\alpha )$ of $[j]\times [k]$ it generates. 
\begin{example}
  Taking $n=30$, $j=5$, $k=4$, and
  $$
  \alpha
  =
\left[ 
\begin{array}{c} 
\alpha_1 
\\ 
\alpha_2 
\\ 
\alpha_3 
\\ 
\alpha_4 
\\ 
\alpha_5 
\\ 
\end{array} 
\right] 
  =
\left[ 
\begin{array}{cccc}
  26 & 15 & 25 & 23
  \\
19 & 23 & 17 & 5\\
24 & 18 & 12 & 20\\
15 & 17 & 7 & 2\\
2 & 26 & 27 & 30 
\end{array} 
\right] 
,
$$ 
the partition $\sqcap (\alpha )$ of
$[5]\times [4]$ is given
in Figure~\ref{fig:diagram1}
 by %
  \begin{align*}
    \sqcap (\alpha )=\big\{ &\{(1,1),(5,2)\},
    \{(1,2),(4,1)\},
    \{(1,3)\},
    \{(1,4),(2,2)\},
    \{(2,1)\},
    \{(2,3),(4,2)\}, 
    \\
    &\{(2,4)\},
    \{(3,1)\}, 
    \{(3,2)\}, 
    \{(3,3)\}, 
    \{(3,4)\}, 
    \{(4,3)\},
    \{(4,4),(5,1)\},
    \{(5,3)\},
    \{(5,4)\}
    \big\}. 
  \end{align*}
\end{example}

\begin{figure}[H]
  \captionsetup[subfigure]{font=footnotesize}
  \centering
  \subcaptionbox{$\alpha =[\alpha_1,\dots,\alpha_5]^\top$.}[.5\textwidth]{%
  \begin{tikzpicture}[scale=0.9] 
  \draw[black, thick] (0,0) rectangle (5,6);

  \node[anchor=east,font=\footnotesize] at (0.6,5) {1};
  \node[anchor=east,font=\footnotesize] at (0.6,4) {2};
  \node[anchor=east,font=\footnotesize] at (0.6,3) {3};
  \node[anchor=east,font=\footnotesize] at (0.6,2) {4};
  \node[anchor=east,font=\footnotesize] at (0.6,1) {5};
  
  \node[anchor=south,font=\footnotesize] at (1,0) {1};
  \node[anchor=south,font=\footnotesize] at (2,0) {2};
  \node[anchor=south,font=\footnotesize] at (3,0) {3};
  \node[anchor=south,font=\footnotesize] at (4,0) {4};
  
  \foreach \i in {1,...,5}
         {
\draw[thick,black] (1,\i) circle (10pt);
\draw[thick,black] (2,\i) circle (10pt);
\draw[thick,black] (3,\i) circle (10pt);
\draw[thick,black] (4,\i) circle (10pt);
} 
\node[blue,font=\normalsize] at (1,1) {2};
\node[blue,font=\normalsize] at (2,1) {26};
\node[blue,font=\normalsize] at (3,1) {27};
\node[blue,font=\normalsize] at (4,1) {30};
\node[blue,font=\normalsize] at (1,2) {15};
\node[blue,font=\normalsize] at (2,2) {17};
\node[blue,font=\normalsize] at (3,2) {7};
\node[blue,font=\normalsize] at (4,2) {2};
\node[blue,font=\normalsize] at (1,3) {24};
\node[blue,font=\normalsize] at (2,3) {18};
\node[blue,font=\normalsize] at (3,3) {12};
\node[blue,font=\normalsize] at (4,3) {20};
\node[blue,font=\normalsize] at (1,4) {19};
\node[blue,font=\normalsize] at (2,4) {23};
\node[blue,font=\normalsize] at (3,4) {17};
\node[blue,font=\normalsize] at (4,4) {5};
\node[blue,font=\normalsize] at (1,5) {26};
\node[blue,font=\normalsize] at (2,5) {15};
\node[blue,font=\normalsize] at (3,5) {25};
\node[blue,font=\normalsize] at (4,5) {23};

 \begin{pgfonlayer}{background}
    \filldraw [line width=4mm,black!3]
      (0.2,0.2)  rectangle (4.8,5.8);
  \end{pgfonlayer}
  \end{tikzpicture}}%
  \subcaptionbox{Partition $\sqcap (\alpha )$.}[.5\textwidth]{
  \begin{tikzpicture}[scale=0.9] 
  \draw[black, thick] (0,0) rectangle (5,6);
  
  \node[anchor=east,font=\footnotesize] at (0.7,5) {1};
  \node[anchor=east,font=\footnotesize] at (0.7,4) {2};
  \node[anchor=east,font=\footnotesize] at (0.7,3) {3};
  \node[anchor=east,font=\footnotesize] at (0.7,2) {4};
  \node[anchor=east,font=\footnotesize] at (0.7,1) {5};
  
  \node[anchor=south,font=\footnotesize] at (1,0) {1};
  \node[anchor=south,font=\footnotesize] at (2,0) {2};
  \node[anchor=south,font=\footnotesize] at (3,0) {3};
  \node[anchor=south,font=\footnotesize] at (4,0) {4};
  
  \filldraw [black] (1,1) circle (2pt);
  \filldraw [black] (2,1) circle (2pt);
  \filldraw [black] (3,1) circle (2pt);
  \filldraw [black] (4,1) circle (2pt);
  \filldraw [black] (1,2) circle (2pt);
  \filldraw [black] (2,2) circle (2pt);
  \filldraw [black] (3,2) circle (2pt);
  \filldraw [black] (4,2) circle (2pt);
  \filldraw [black] (1,3) circle (2pt);
  \filldraw [black] (2,3) circle (2pt);
  \filldraw [black] (3,3) circle (2pt);
  \filldraw [black] (4,3) circle (2pt);
  \filldraw [black] (2,3) circle (2pt);
  \filldraw [black] (1,4) circle (2pt);
  \filldraw [black] (2,4) circle (2pt);
  \filldraw [black] (3,4) circle (2pt);
  \filldraw [black] (4,4) circle (2pt);
  \filldraw [black] (1,5) circle (2pt);
  \filldraw [black] (2,5) circle (2pt);
  \filldraw [black] (3,5) circle (2pt);
  \filldraw [black] (4,5) circle (2pt);
  
  \draw[very thick] (1,2) -- (2,5);
  \draw[very thick] (1,1) -- (4,2);
  \draw[very thick] (1,5) -- (2,1);
  \draw[very thick] (2,2) -- (3,4);
  \draw[very thick] (2,4) -- (4,5);
 \begin{pgfonlayer}{background}
    \filldraw [line width=4mm,black!3]
      (0.2,0.2)  rectangle (4.8,5.8);
  \end{pgfonlayer}
  \end{tikzpicture}}%
  
  \caption{Example for the mapping $\sqcap$ with $j=5$ and $k=4$.}
  \label{fig:diagram1}
  \end{figure}
  
  \vspace{-0.7cm}

  \noindent
   We let $v(G):=|V_G|$ and $e(G):=|E_G|$ be the number of vertices and the number of edges of any graph 
 $G=(V_G,E_G)$ with vertex set $V_G$ and edge set $E_G$.
A subgraph of $G$ is a graph $H=(V_{H},E_{H})$ such that $V_H\subset V_G$ and $E_H\subset E_G$, and $H$ is an induced subgraph of $G$, if $E_H$ consists of all edges of $G$ having both endpoints in $V_H$.
Two graphs $G=(V_G,E_G)$ and $H=(V_H,E_H)$ are isomorphic
if there is a bijection $T:V_G\to V_H$ such that $\{u,v\}\in E_G$ if and only if $\{T(u),T(v)\}\in E_H$ for any $u\ne v\in V_G$, in which case we write $H=T(G)$.
The permutations $\alpha$ of $V_G$ such that $\alpha(G)=G$ form a group called the automorphism group,
and we let $a(G)$ denote the cardinality of this group.

\begin{definition}\label{def:diagram-graph}
  Consider $G_1,\ldots , G_j$ copies of
  a connected graph $G$ with $v(G) = k\geq 2$ vertices,   
  respectively built on
  $\pi_1,\ldots , \pi_j$, $j\geq 1$,
  and let $\rho \in\Pi( [j] \times[k])$
  be a partition of $ [j] \times [k]$. 
\begin{enumerate}[\rm 1.]
\item
  We let $\widebar{\rho}_G$ denote the
  contraction multigraph of the graph $G^{\otimes j}$ 
  constructed on the blocks of $\rho$ 
  by adding an edge between two blocks $\rho_1,\rho_2$ of the
  partition $\rho$ whenever there exist $(i,l_1)\in \rho_1$
  and $(i,l_2)\in \rho_2$ such that $(l_1,l_2)$ is an edge in $G_i$.
\item 
  We let $\rho_G$ be the graph constructed
  on the blocks of $\rho$
  by removing redundant edges
  in $\widebar{\rho}_G$,
  so that at most one edge remains between any two blocks $\rho_1,\rho_2\in\rho$. 
 \end{enumerate}
\end{definition}

\begin{example}
  \label{fjkld243}
   Consider $\rho\in\Pi([3]\times[3])$ of the form 
\begin{align*} 
  \rho & =
  \big\{
  \{(1,1),(2,1)\},
  \{(1,2),(3,2)\},
  \{(1,3),(3,3)\},
  \{(2,3)\},
  \{(2,2),(3,1)\}  
  \big\}
  \\
  &
  =
  \{
  \rho_1,
  \rho_2,
  \rho_3,
  \rho_4,
  \rho_5
  \}, 
\end{align*} 
  and let $G$ be a triangle on the vertex set $V_G=[3]$. 
  Figure~\ref{fig:diagram1-3} presents the multigraph $\widebar{\rho}_G$
  and corresponding graph $\rho_G$, 
  where $\rho \in\Pi([3]\times[3])$
  is non-flat and connected. 
  
\smallskip

\begin{figure}[H]
\captionsetup[subfigure]{font=footnotesize}
\centering
\subcaptionbox{Multigraph $\widebar{\rho}_G$ before merging edges and vertices.}[.5\textwidth]{%
\begin{tikzpicture}%
\draw[step=1cm, very thin, gray!40] (1,0) grid (5,4);
\draw[black, thick] (1,0) rectangle (5,4);
\foreach \i in {2,3}
{
\filldraw [black] (2,\i) circle (2pt);
\filldraw [black] (3,\i) circle (2pt);
\filldraw [black] (4,\i) circle (2pt);
\draw[thick, dash dot,blue] (2,\i) -- (3,\i);
\draw[thick, dash dot,blue] (3,\i) -- (4,\i);
\draw[thick, dash dot,blue] (2,\i) .. controls (3,\i-0.5) .. (4,\i);
}

\filldraw [black] (2,1) circle (2pt);
\filldraw [black] (3,1) circle (2pt);
\filldraw [black] (4,1) circle (2pt);
\draw[thick, dash dot,blue] (2,1) -- (3,1);
\draw[thick, dash dot,blue] (3,1) -- (4,1);
\draw[thick, dash dot,blue] (2,1) .. controls (3,1+0.5) .. (4,1);

\node[anchor=north,font=\scriptsize] at (2,1) {(3,1)};
\node[anchor=north,font=\scriptsize] at (3,1) {(3,2)};
\node[anchor=north,font=\scriptsize] at (4,1) {(3,3)};
\node[anchor=south,font=\scriptsize] at (2,3) {(1,1)};
\node[anchor=south,font=\scriptsize] at (3,3) {(1,2)};
\node[anchor=south,font=\scriptsize] at (4,3) {(1,3)};
\node[anchor=south,font=\scriptsize] at (1.5,1.7) {(2,1)};
\node[anchor=south,font=\scriptsize] at (3,2) {(2,2)};
\node[anchor=south,font=\scriptsize] at (4.5,1.7) {(2,3)};

\draw[thick] (3,1) .. controls (3.5,2) .. (3,3);
\draw[thick] (2,3) -- (2,2);
\draw[thick] (3,2) -- (2,1);
\draw[thick] (4,1) .. controls (3.5,2) .. (4,3);
\begin{pgfonlayer}{background}
    \filldraw [line width=4mm,black!3]
      (1.2,0.2)  rectangle (4.8,3.8);
  \end{pgfonlayer}
\end{tikzpicture}}%
\subcaptionbox{Graph $\rho_G$ after merging edges and vertices.}[.5\textwidth]{
\begin{tikzpicture}%
\draw[step=1cm, very thin, gray!40] (1,0) grid (5,4);
  \draw[black, thick] (1,0) rectangle (5,4);
  \filldraw [black] (2,2) circle (2pt);
  \filldraw [black] (3,1) circle (2pt);
  \filldraw [black] (3,3) circle (2pt);
  \filldraw [black] (4,1) circle (2pt);
  \filldraw [black] (4,2) circle (2pt);
  \filldraw [black] (4,1) circle (2pt);
  \draw[thick,blue] (2,2) -- (3,3);
  \draw[thick,blue] (3,3) -- (3,1);
  \draw[thick,blue] (2,2) -- (3,1);
  \draw[thick,blue] (3,3) -- (4,2);
  \draw[thick,blue] (2,2) -- (4,2);
  \draw[thick,blue] (3,1) -- (4,2);
  \draw[thick,blue] (3,1) -- (4,1);
  \draw[thick,blue] (2,2) -- (4,1) -- (4,1);

  \node[anchor=south,font=\footnotesize] at (3,3) {$\rho_2$};
  \node[anchor=north,font=\footnotesize] at (3,1) {$\rho_5$};
  \node[anchor=east,font=\footnotesize] at (2,2) {$\rho_1$};
  \node[anchor=west,font=\footnotesize] at (4,2) {$\rho_3$};
  \node[anchor=north,font=\footnotesize] at (4,1) {$\rho_4$};

  \begin{pgfonlayer}{background}
    \filldraw [line width=4mm,black!3]
      (1.2,0.2)  rectangle (4.8,3.8);
  \end{pgfonlayer}
\end{tikzpicture}}%
\caption{
 Example of graph $\rho_G$ with $j=3$ and $k=3$.}
\label{fig:diagram1-3}
\end{figure}

\vspace{-0.7cm}

\end{example}

\noindent
  The contraction multigraph $\widebar{\rho}_G$
  constructed in Definition~\ref{def:diagram-graph}
  is also denoted by $G^{\otimes j}$
  in \cite[Page~104]{feray16}, where $G^{\otimes j}$ stands for
  the graph made of $j$ components isomorphic to~$G$.

We will use the following standard notation for the asymptotic behavior of the relative order of magnitude of two functions $f(n)$ and $g(n )>0$ as $n$ tends to infinity.
We write
\begin{itemize}
\item
$f(n ) \gtrsim g(n )$ if $\liminf_{n \to\infty} f(n ) / g(n )>0$,
\item $f(n)\ll g(n)$, or $g(n)\gg f(n)$,
   if $f(n )\geq 0$, $g(n)>0$, and $\lim_{n\to \infty} f(n)/g(n) = 0$. 
\end{itemize}
\noindent 
We close this section by recalling the following lemma, see Lemma~2.8 in \cite{LiuPrivault}
or Proposition 6.1 in \cite{schulte-thaele}. 
       \begin{lemma}
         \label{fkla12}
          \noindent
        $a)$ The cardinality of the set \ ${\rm CNF}(j,k)$
       of connected non-flat partitions of $[j]\times[k]$ satisfies 
       \begin{equation}
         \nonumber %
        |{\rm CNF}(j,k)| \leq j!^k k!^{j-1}, 
        \qquad j,k \geq 1. 
      \end{equation}
      \noindent
      $b)$ 
       The cardinality of the set 
      $ \mathcal{M}(j,k)$
       of maximal connected non-flat partitions of $[j]\times[k]$
       satisfies 
      \begin{equation*}
        |\mathcal{M}(j,k)|=k^{j-1}\prod_{i=1}^{j-1}(1+(k-1)i),
        \qquad j,k\geq 1, 
      \end{equation*}
       with the bounds 
      \begin{equation*}
          ( (k-1)k )^{j-1}(j-1)!\le
          |\mathcal{M}(j,k)|
           \leq ( (k-1)k )^{j-1}j!, \quad j\geq 1, \ k\geq 2. 
      \end{equation*}
      \end{lemma} 
\section{Moments of generalized $U$-statistics}
\label{s3}
 \noindent
In this section, we prove a moment identity
for generalized $U$-statistics, see \cite{Janson91} and \cite[Chapter 11]{janson},
 by combining multiple integrals with set partitions.
 Let $\mathcal{M}$ and $\mathcal{S}$ 
 be respectively a measurable space and a Borel space. 
\begin{definition} 
  Given $f:\mathcal{S}^k\times \mathcal{M}^{(k-1)k/2}\to\R$ a bounded measurable
  function, $k\geq 2$,
 $\mathcal{Y}=\{Y_{i ,j }\}_{1\leq i < j \leq n}$
 a sequence of independent $\mathcal{M}$-valued 
 random variables with common probability distribution $\mathbb{Q}$, and
 $\mathcal{X}=\{X_i \}_{1\leq i \leq n}$
 a sequence of independent $\mathcal{S}$-valued 
 random variables with continuous distribution $\mu$, we let
 $S_{n,k}(f)$ denote the generalized $U$-statistics
 defined as 
\begin{equation}
\label{def:gUstat}
  S_{n,k}(f):=\sum_{\beta\in[n]^k_\ne }f (X_{\beta(1)},\dots,X_{\beta(k)},Y_{\beta(1),\beta(2)},\dots,Y_{\beta(1),\beta(k)},Y_{\beta(2),\beta(3)},\dots,Y_{\beta(k-1),\beta(k)} ), 
\end{equation}
where $Y_{i ,j }=Y_{j ,i }$ if $j >i$.
\end{definition} 
Note that the function $f$ is not
required to be symmetric as in \cite{Janson91}
or \cite{zhangzs}.
 For any $\beta =(\beta (1),\dots,\beta (k))\in[n]^k_{\ne}$, 
 we also let 
 $$
 f(\mathbf{X}_\beta ,\mathbf{Y}_\beta )
 = f\big(X_{\beta (1)},\dots,X_{\beta (k)},Y_{\beta (1),\beta (2)},\dots,Y_{\beta (k-1),\beta (k)}\big)
 $$
 for shortness of notation. 
 Letting $\rho_{K_k}$ denote the contraction graph
 of the complete graph $K_k$ on $k$ vertices,
 see Definition~\ref{def:diagram-graph},
 for any set partition
$\rho =\{b_1,\dots,b_{|\rho |}\}\in\Pi([j]\times[k])$ we 
consider the function 
$$
\left(\bigotimes_{i=1}^j f
\right)_{\! \! \! \rho }
 \ \! : \ \! \mathcal{S}^{|\rho |}\times \mathcal{M}^{e(\rho_{K_k})}\to\R
$$
 defined as 
\begin{eqnarray*}
\left(\bigotimes_{i=1}^j f
\right)_{\! \! \! \rho }
\big(
\{x_\nu\}_{1\leq \nu\leq |\rho |},\{y_{(\nu,\upsilon)}\}_{(\nu,\upsilon)\in E_{\rho_{K_k}}}
\big) :=\prod_{i=1}^jf\big(x^{(i)}_1,\dots,x^{(i)}_k ,y^{(i)}_{1,2},\dots,y^{(i)}_{k-1,k}\big),
\end{eqnarray*}
 where $x^{(i)}_\ell:=x_\nu$ if $(i,\ell)\in b_\nu$ for
 $1\leq \nu\leq |\rho |$, 
 and $y^{(i)}_{u,v}:=y_{(t,s)}$ if $(i,u)\in b_t$ and $(i,v)\in b_s$,
 $s,t=1,\ldots , |b|$.
 In Theorem~\ref{jklfd12} we 
 provide a moment identity for the generalized $U$-statistics $S_{n,k}(f)$
 using the partition diagram
 language of Definition~\ref{def:diagram-graph}.
  \begin{thm}
    \label{jklfd12} 
  Let $f:\mathcal{S}^k\times \mathcal{M}^{(k-1)k/2}\to\R$ be a bounded measurable function with $k\geq 2$.
  For any $j,n\ge1$, we have 
  \begin{equation}
    \nonumber %
    \E\left[(S_{n,k}(f))^j\right]=\sum_{\substack{\rho\in\Pi([j]\times[k])\\\rho\wedge\pi=\widehat{0}}}
           \frac{n!}{(n-|\rho|)!}
           \int_{\mathcal{S}^{|\rho|}\times\mathcal{M}^{e(\rho_{K_k})}}\left(\bigotimes_{i=1}^j f\right)_{\! \! \! \rho} (\mathbf{x},\mathbf{y})\mu^{\otimes |\rho|}
           ( \mathrm{d}\mathbf{x} )
           \mathbb{Q}^{e(\rho_{K_k})}
            ( \mathrm{d}\mathbf{y} ) .
  \end{equation}
\end{thm}
\begin{Proof}
 From the definition of
    $\sqcap( \alpha )$, we have 
\begin{align*}
  \E\left[(S_{n,k}(f))^j\right]&= \sum_{\alpha_1\in [n]^k_{\ne},\ldots,\alpha_j\in[n]^k_{\ne}}\E\left[\prod_{i=1}^jf(\mathbf{X}_{\alpha_i},\mathbf{Y}_{\alpha_i})\right]
    \\
    &
      = \sum_{\substack{\rho\in\Pi([j]\times[k])\\\rho\wedge\pi=\widehat{0}}}
    \sum_{\alpha_1\in [n]^k_{\ne},\ldots,\alpha_j\in[n]^k_{\ne}
    \atop \sqcap (\alpha ) = \rho }
    \E\left[\prod_{i=1}^jf(\mathbf{X}_{\alpha_i},\mathbf{Y}_{\alpha_i})\right]
       \\
    &
      = \sum_{\substack{\rho\in\Pi([j]\times[k])\\\rho\wedge\pi=\widehat{0}}}
    \sum_{\alpha_1\in [n]^k_{\ne},\ldots,\alpha_j\in[n]^k_{\ne}
    \atop \sqcap (\alpha ) = \rho }
    \int_{\mathcal{S}^{|\rho|}\times\mathcal{M}^{e(\rho_{K_k})}}\left(\bigotimes_{i=1}^j f\right)_{\! \! \! \rho} (\mathbf{x},\mathbf{y})\mu^{\otimes |\rho|}
  ( \mathrm{d}\mathbf{x} )
  \mathbb{Q}^{e(\rho_{K_k})} (
  \mathrm{d}\mathbf{y} )
\\
  &= \sum_{\substack{\rho\in\Pi([j]\times[k])\\\rho\wedge\pi=\widehat{0}}}C_n(\rho)\int_{\mathcal{S}^{|\rho|}\times\mathcal{M}^{e(\rho_{K_k})}}\left(\bigotimes_{i=1}^j f\right)_{\! \! \! \rho} (\mathbf{x},\mathbf{y})\mu^{\otimes |\rho|}
  ( \mathrm{d}\mathbf{x} )
  \mathbb{Q}^{e(\rho_{K_k})} (
  \mathrm{d}\mathbf{y} )
,
\end{align*}
  where %
\begin{equation}
\nonumber %
  C_n(\rho):=
  \frac{n!}{(n-|\rho|)!},
  \quad
  \rho\in\Pi([j] \times [k]), 
\end{equation}
 denotes the count of the 
 $\alpha  = [\alpha_1,\ldots , \alpha_j]^\top \in [n]^k_{\ne}\times \cdots \times [n]^k_{\ne}$ such that $\sqcap (\alpha )=\rho$,
 $j,k\ge1$, with $\rho$ non-flat. 
\end{Proof}
In particular, since $\widehat{0}\in\Pi([1]\times[k])$ is the only
non-flat partition of $[1]\times[k]$, %
 we have %
 \begin{equation}
   \label{exp-1}
    \E\left[S_{n,k}(f)\right]=
    \frac{n!}{(n-k)!}
    \int_{\mathcal{S}^k\times\mathcal{M}^{(k-1)k/2}}f(\mathbf{x},\mathbf{y})\mu^{\otimes k}
    ( \mathrm{d}\mathbf{x})
    \mathbb{Q}^{(k-1)k/2}
    ( \mathrm{d}\mathbf{y} ). 
\end{equation}
\subsubsection*{Moments of standard $U$-statistics}
\noindent
 Given $j\ge1$ and $f^{(i)}:(\R^d)^k\to\R$, $i=1,\dots,j$,
       measurable functions, we let 
       \begin{equation}
         \nonumber %
       \left(\bigotimes_{i=1}^jf^{(i)} \right)(x_{1,1},\ldots,x_{1,k},\ldots,
       x_{j,1},\ldots , x_{j,k}):=\prod_{i=1}^jf^{(i)} (x_{i,1},\dots,x_{i,k}), 
\end{equation}
  and for $\rho\in \Pi([j] \times [k] )$ we denote by 
        $$
        \left(
        \bigotimes_{i=1}^jf^{(i)} \right)_{\! \! \! \rho} \ \! : \ \! (\R^d)^{|\rho|}\to\R
        $$
        the function obtained by equating any two variables
        whose indexes belong to a same block of $\rho$.
        In Corollary~\ref{momid-1},
        as a consequence of Theorem~\ref{jklfd12}, 
        we obtain a moment identity for
        standard $U$-statistics of order $k\geq 1$,
 of the form 
\begin{equation}
  \nonumber %
  S_n(f):=\sum_{\beta \in[n]^k_{\ne}}f\big(X_{\beta (1)},\dots,X_{\beta (k)}\big). 
\end{equation}
\begin{corollary}\label{momid-1}
  Let $f:(\R^d)^k\to\R$, $k\ge1$, be a bounded,
   not necessarily symmetric,
 measurable function. 
    For any $j,n\ge1$, we have 
  \begin{eqnarray}
    \E\left[ ( S_n(f) )^j\right]=\sum_{\substack{\rho\in\Pi([j]\times[k])\\\rho\wedge\pi=\widehat{0}}}
    \frac{n!}{(n-|\rho|)!}
    \int_{(\R^d)^{|\rho|}}\left(\bigotimes_{i=1}^j f\right)_{\! \! \! \rho}
    ( x_1, \ldots , x_{|\rho|} )
    \mu ( \mathrm{d}x_1 ) \cdots\mu ( \mathrm{d}x_{|\rho|} ). %
\nonumber %
  \end{eqnarray}
\end{corollary}
\begin{Proof}
  Since the sequence $(X_1,\dots,X_n)$ is i.i.d., we have
  \begin{align*}
    \E\left[(S_n(f))^j\right] %
    &=\sum_{\alpha_1 \in [n]^k_{\ne},\ldots,\alpha_j \in [n]^k_{\ne}}\E\left[\prod_{i=1}^jf\big(X_{\alpha_i(1)},\dots,X_{\alpha_i(k)}\big)\right]\\
    &= \sum_{\substack{\rho\in\Pi([j]\times[k])\\\rho\wedge\pi=\widehat{0}}}C_n(\rho)\int_{(\R^d)^{|\rho|}}\left(\bigotimes_{i=1}^j f\right)_{\! \! \! \rho}
    ( x_1, \ldots , x_{|\rho|} )
    \mu ( \mathrm{d}x_1 ) \cdots\mu ( \mathrm{d}x_{|\rho|} ). 
\end{align*}
\end{Proof}
\begin{remark}
  For the $V$-statistics 
  \begin{equation}
    \nonumber %
    V_{n,k}(f):=\sum_{\beta \in[n]^k}f\big(X_{\beta (1)},\dots,X_{\beta (k)}\big) 
  \end{equation}
  with possible repeated indices,
  for $j,n\ge1$, we have the similar moment identity 
  \begin{equation}
    \label{moment-2}
    \E\left[(V_{n,k}(f))^j\right]=\sum_{\rho\in\Pi([j]\times[k])}
    \frac{n!}{(n-|\rho|)!}
    \int_{(\R^d)^{|\rho|}}\left(\bigotimes_{i=1}^j f\right)_{\! \! \! \rho}
    ( x_1, \ldots , x_{|\rho|} )
    \mu ( \mathrm{d}x_1 ) \cdots\mu ( \mathrm{d}x_{|\rho|} ).
  \end{equation}
  The proof of \eqref{moment-2} is almost the same, except for
  the removal of the non-flat restriction, i.e.  
  \begin{align*}
    \E\left[(V_{n,k}(f))^j\right]&= \sum_{\alpha_1,\dots,\alpha_j\in[n]^k}\E\left[\prod_{i=1}^jf\big(X_{\alpha_i(1)},\dots,X_{\alpha_i(k)}\big)\right]
    \\
    &= \sum_{\rho\in\Pi([j]\times[k])}
    \frac{n!}{(n-|\rho|)!}
    \int_{(\R^d)^{|\rho|}}
    \left(\bigotimes_{i=1}^j f
    \right)_{\! \! \! \rho}
    ( x_1, \ldots , x_{|\rho|} )
    \mu ( \mathrm{d}x_1 ) \cdots\mu ( \mathrm{d}x_{|\rho|} ).
  \end{align*}
\end{remark}
\section{Cumulant bounds for generalized $U$-statistics}
\label{s4}
\noindent 
To study the asymptotic behaviour of $U$-statistics and generalized $U$-statistics, a common practice is to use Hoeffding decompositions \cite{hoeffding61}, \cite{mandelbaum}. Through the orthogonal decomposition, one finds the asymptotic distribution of $U$-statistics is determined by the ``smallest'' component
appearing in its decomposition, see \cite[Lemma~2]{Janson91} and also \cite[Theorem~11.3]{janson}.
 In what follows, we consider the case where the asymptotic distribution of $S_{n,k}(f)$ is normal.
    Assumption~\ref{assu1}
   will be needed for the derivation
   of cumulant estimates for the generalized $U$-statistics $S_{n,k}(f)$
   in Theorem~\ref{jkldd12}.
 \begin{assumption}\label{assu1}
  Let $f:\mathcal{S}^k\times \mathcal{M}^{(k-1)k/2}\to\R$ be a bounded measurable function with $k\geq 2$. We assume that 
  \begin{equation}
 \label{assu01}
       \Var
      \left[\sum_{\ell=1}^k
        f_{(\ell)} (X_1) \right] >0. 
  \end{equation}
  where 
  \begin{eqnarray}
    \nonumber %
    f_{(i)}(x):=\int_{\mathcal{S}^{k-1}\times \mathcal{M}^{(k-1)k/2}}f(\mathbf{x},\mathbf{y})\mu^{\otimes(k-1)}\{\mathrm{d}x_1\cdots\mathrm{d}x_{i-1}\mathrm{d}x_{i+1}\cdots\mathrm{d}x_k\}\mathbb{Q}^{(k-1)k/2}\{\mathrm{d}\mathbf{y}\}
 \end{eqnarray}
for $i=1,\dots,k$,
and $\mathbf{x}=(x_1,\dots,x_{i-1},x,x_{i+1},\dots,x_k )$.
\end{assumption}
 Assumption~\ref{assu1} amounts to saying
 that the principal degree of $f$ equals $1$, see \cite{Janson91}.
In most of the existing literature, including \cite{Janson91}, \cite{KaurRollin21}, 
\cite{zhangzs}, the asymptotic normality of generalized $U$-statistics relies on orthogonal decomposition on certain $L_2$ spaces. However, in practice the orthogonal decomposition of count statistics can be
 intractable in general. 
\begin{thm}
    \label{jkldd12} 
  Let $f:\mathcal{S}^k\times \mathcal{M}^{(k-1)k/2}\to\R$ be a measurable function, with $k\geq 2$.
  Then, for $n\ge1$
  the $j$-$th$ cumulant of the generalized $U$-statistics $S_{n,k}(f)$
  satisfies the bound   
  \begin{equation}
    \label{fjkldf}
     | \kappa_j(S_{n,k}(f))|\leq n^{1+(k-1)j}
    \Vert f\Vert_\infty^j j^{j-1} (j!)^k (k!)^{j-1}, \qquad j\geq 1.
  \end{equation}
  In addition, if the function 
  $f$ satisfies Assumption~\ref{assu1}, then we have 
  \begin{equation}
    \label{lowb-2}
   \kappa_2(S_{n,k}(f))\geq \frac{C n!}{(n-2k+1)!},
   \quad
   n \geq N(f,k),
 \end{equation}
   where
   $C(f,k)>0$ and 
   $N(f,k) \geq 1$ depend only on $f$ and on $k \geq 2$.
\end{thm}
\begin{Proof}
  From Property~C in \cite[Page~29]{MalyshevMinlos91}, we know that
  if $\alpha = [\alpha_1,\ldots , \alpha_j]^\top \in [n]^k_{\ne}\times \cdots \times [n]^k_{\ne}$ is disconnected, i.e. 
  if there exists a partition $\{A,B\}$ of $[j]$ such that 
  \begin{equation*}
    \big\{\alpha_i(\ell)\}_{(i, \ell ) \in A \times [k]}
    \cap
     \big\{\alpha_i(\ell)\}_{(i,\ell ) \in B \times [k] }=\emptyset,
  \end{equation*}
  then the joint cumulant $\kappa(f(\mathbf{X}_{\alpha_1},\mathbf{Y}_{\alpha_1}),\dots,f(\mathbf{X}_{\alpha_j},\mathbf{Y}_{\alpha_j}))$
  vanishes.
   Hence, we have 
  \begin{align}
    \nonumber
    \kappa_j(S_{n,k}(f))&= \sum_{\alpha_1,\dots,\alpha_j\in[n]^k_{\ne}}\kappa(f(\mathbf{X}_{\alpha_1},\mathbf{Y}_{\alpha_1}),\dots,f(\mathbf{X}_{\alpha_j},\mathbf{Y}_{\alpha_j}))
    \\
    \nonumber
    &= \sum_{\substack{\alpha_1  \in [n]^k_{\ne},\dots,\alpha_j \in [n]^k_{\ne}
        \\\mathrm{connected}}}\kappa(f(\mathbf{X}_{\alpha_1},\mathbf{Y}_{\alpha_1}),\dots,f(\mathbf{X}_{\alpha_j},\mathbf{Y}_{\alpha_j}))
    \\
    \label{fjkldf1} 
        &= \sum_{\rho\in{\rm CNF}(j,k)}\sum_{\substack{
        \alpha = [\alpha_1,\ldots , \alpha_j]^\top \in
        ([n]^k_{\ne})^j
\\\sqcap (\alpha )=\rho}}\kappa(f(\mathbf{X}_{\alpha_1},\mathbf{Y}_{\alpha_1}),\dots,f(\mathbf{X}_{\alpha_j},\mathbf{Y}_{\alpha_j})). \qquad 
      \end{align} 
  By the cumulant-moment relation \eqref{cummom1} we have,
  for any $[\alpha_1,\ldots , \alpha_j]^\top \in
  [n]^k_{\ne} \times \cdots \times [n]^k_{\ne}$, 
\begin{align}
  \nonumber
  \big|\kappa \big(
  f(\mathbf{X}_{\alpha_1},\mathbf{Y}_{\alpha_1}),\dots,f(\mathbf{X}_{\alpha_j},\mathbf{Y}_{\alpha_j})
  \big)\big|
  &\leq  \sum_{\sigma=\{b_1,\dots,b_l\}\in\Pi([j])}(l-1)!\prod_{i=1}^l\left|\E\left[\prod_{\ell\in b_i}f(\mathbf{X}_{\alpha_\ell},\mathbf{Y}_{\alpha_\ell})\right]\right|
  \\
  \nonumber  &\leq \sum_{\sigma=\{b_1,\dots,b_l\}\in\Pi([j])}(l-1)! \Vert f\Vert_\infty^j
  \\
  \nonumber
  &=\Vert f\Vert_\infty^j\sum_{l=1}^j(l-1)!S(j,l)
  \\
  \nonumber
  &\leq \Vert f\Vert_\infty^j \frac{1}{j} \sum_{l=1}^j
  \frac{j!}{(j-l)!}
  S(j,l)
  \\
\label{uppbou1}
  &= \Vert f\Vert_\infty^j j^{j-1},
\end{align}
where $S(j,l)$ is the Stirling number of the second kind. 
  Therefore, from \eqref{fjkldf1} we obtain 
  \begin{align*}
    \kappa_j(S_{n,k}(f))& \leq 
    \sum_{\rho\in{\rm CNF}(j,k)}
    \frac{n!}{(n-|\rho|)!}
    \Vert f\Vert_\infty^j j^{j-1}
    \\ 
    &=\sum_{r=k}^{1+(k-1)j}
    \frac{n!}{(n-r)!}
    \sum_{\substack{\rho\in{\rm CNF}(j,k)\\|\rho|=r}}
        \Vert f\Vert_\infty^j j^{j-1}
      \\
  & \leq n^{1+(k-1)j}\Vert f\Vert_\infty^j j^{j-1}|{\rm CNF}(j,k)|, 
 \end{align*}
 which yields \eqref{fjkldf} from Lemma~\ref{fkla12}-$(a)$. %
 Letting $j = 2$ in \eqref{fjkldf1}, we have 
  \begin{align} 
    \nonumber
    \kappa_2(S_{n,k}(f))&= \sum_{\rho\in{\rm CNF}(2,k)}\sum_{\substack{
      \alpha = [\alpha_1,\alpha_2]^\top \in ([n]^k_{\ne})^2
      \\\sqcap (\alpha )=\rho}}\kappa(f(\mathbf{X}_{\alpha_1},\mathbf{Y}_{\alpha_1}),f(\mathbf{X}_{\alpha_2},\mathbf{Y}_{\alpha_2}))
    \\
    \nonumber
    &= \sum_{\rho\in{\rm CNF}(2,k)}\sum_{\substack{
      \alpha = [\alpha_1,\alpha_2]^\top \in ([n]^k_{\ne})^2
      \\\sqcap (\alpha )=\rho}}\cov(f(\mathbf{X}_{\alpha_1},\mathbf{Y}_{\alpha_1}),f(\mathbf{X}_{\alpha_2},\mathbf{Y}_{\alpha_2}))
    \\
      \label{decp1}
      &=
      R +
      \sum_{\substack{\rho\in{\rm CNF}(2,k)\\|\rho|=2k-1}}\sum_{\substack{
      \alpha = [\alpha_1,\alpha_2]^\top \in ([n]^k_{\ne})^2\\
      \sqcap (\alpha )=\rho}}\cov(f(\mathbf{X}_{\alpha_1},\mathbf{Y}_{\alpha_1}),f(\mathbf{X}_{\alpha_2},\mathbf{Y}_{\alpha_2})),
  \end{align} 
  where 
  \begin{equation}
    \nonumber
    R:=\sum_{\substack{\rho\in{\rm CNF}(2,k)\\|\rho|\le2k-2}}\sum_{\substack{
      \alpha = [\alpha_1,\alpha_2]^\top \in ([n]^k_{\ne})^2
     \\\sqcap (\alpha )=\rho}}\cov(f(\mathbf{X}_{\alpha_1},\mathbf{Y}_{\alpha_1}),f(\mathbf{X}_{\alpha_2},\mathbf{Y}_{\alpha_2})).
  \end{equation}
  Because $\{X_1,\dots,X_n\}$ and $\{Y_{i ,j }\}_{1\leq i <j \leq n}$ are
  both i.i.d. random elements, we have, for any
  $\alpha = [\alpha_1,\alpha_2]^\top \in ([n]^k_{\ne})^2$,
  $\beta = [\beta_1,\beta_2]^\top \in ([n]^k_{\ne})^2$
  \begin{equation}
    \nonumber
    \cov(f(\mathbf{X}_{\alpha_1},\mathbf{Y}_{\alpha_1}),f(\mathbf{X}_{\alpha_2},\mathbf{Y}_{\alpha_2}))=\cov(f(\mathbf{X}_{\beta_1},\mathbf{Y}_{\beta_1}),f(\mathbf{X}_{\beta_2},\mathbf{Y}_{\beta_2})),
  \end{equation}
  if $\sqcap (\alpha )=\sqcap (\beta )$. Let $\sigma\in{\rm CNF}(2,k)$ with $|\sigma|=2k-1$ and $\{(1,s),(2,t)\}\in\sigma$, $1\leq s\not= t \leq k$.
  Taking $\alpha = [\alpha_1,\alpha_2]^\top \in ([n]^k_{\ne})^2$
  such that $\sqcap (\alpha )=\sigma$
  we have $\alpha_1(s)=\alpha_2(t)$, hence
\begin{align} 
  \nonumber
  &\cov(f(\mathbf{X}_{\alpha_1},\mathbf{Y}_{\alpha_1}),f(\mathbf{X}_{\alpha_2},\mathbf{Y}_{\alpha_2}))
  \\
  \nonumber
  &\qquad = \E\left[f(\mathbf{X}_{\alpha_1},\mathbf{Y}_{\alpha_1})f(\mathbf{X}_{\alpha_2},\mathbf{Y}_{\alpha_2})\right]
  -
  \E\left[f(\mathbf{X}_{\alpha_1},\mathbf{Y}_{\alpha_1})\right]
  -\E\left[f(\mathbf{X}_{\alpha_2},\mathbf{Y}_{\alpha_2})\right]
  \\
  \nonumber
  &\qquad = \E\left[\E[f(\mathbf{X}_{\alpha_1},\mathbf{Y}_{\alpha_1})f(\mathbf{X}_{\alpha_2},\mathbf{Y}_{\alpha_2})|X_{\alpha_1(t)}]\right]
  -
  \E\left[\E[f(\mathbf{X}_{\alpha_1},\mathbf{Y}_{\alpha_1})|X_{\alpha_1(s)}]\right]
  \E\left[\E[f(\mathbf{X}_{\alpha_2},\mathbf{Y}_{\alpha_2})|X_{\alpha_2(t)}]\right]
  \\
  \nonumber
  &\qquad = \E\left[f_{(s)}(X_{\alpha_1(s)})f_{(t)}(X_{\alpha_1(t)})\right]
  -
  \E\left[f_{(s)}(X_{\alpha_1(s)})\right]
  \E\left[f_{(t)}(X_{\alpha_1(t)})\right]
  \\
  \nonumber
  &\qquad = \E\left[f_{(s)}(X_1)f_{(t)}(X_1)\right]
  -
  \E\left[f_{(s)}(X_1)\right]
  \E\left[f_{(t)}(X_1)\right]
  \\
  \nonumber
  &\qquad
    = \cov (f_{(s)}(X_1), f_{(t)}(X_1)) 
   \\
  \nonumber
  &\qquad =:r_{\sigma}.
\end{align}  
  Therefore, the second term in \eqref{decp1} becomes 
\begin{align} 
  \nonumber
&  \sum_{\substack{\rho\in{\rm CNF}(2,k)\\|\rho|=2k-1}}\sum_{\substack{
      \alpha = [\alpha_1,\alpha_2]^\top \in ([n]^k_{\ne})^2
     \\\sqcap (\alpha )=\rho}}\cov(f(\mathbf{X}_{\alpha_1},\mathbf{Y}_{\alpha_1}),f(\mathbf{X}_{\alpha_2},\mathbf{Y}_{\alpha_2}))
  \\
  \nonumber
  & \qquad \quad = \sum_{\substack{\rho\in{\rm CNF}(2,k)\\|\rho|=2k-1}}\sum_{\substack{
      \alpha = [\alpha_1,\alpha_2]^\top \in ([n]^k_{\ne})^2\\\sqcap (\alpha )=\rho}}r_\rho
  \\
  \nonumber
  & \qquad \quad = \sum_{\substack{\rho\in{\rm CNF}(2,k)\\|\rho|=2k-1}}\frac{n!}{(n-2k+1)!}r_\rho
  \\
  \nonumber
  & \qquad \quad = \frac{n!}{(n-2k+1)!}\sum_{\substack{\rho\in{\rm CNF}(2,k)\\|\rho|=2k-1}}r_\rho
  \\
  \nonumber
  & \qquad \quad
      = \frac{n!}{(n-2k+1)!}
  \sum_{1\leq s,t\leq k} \cov (f_{(s)}(X_1), f_{(t)}(X_1)) 
    \\
\label{decp2}
& \qquad \quad = \frac{n!}{(n-2k+1)!}\Var\left[
  \sum_{\ell=1}^kf_{(\ell)}(X_1)\right].
\end{align} 
 On the other hand, from \eqref{uppbou1}, we have 
\begin{align} 
  \nonumber
  |R|&\leq \sum_{\substack{\rho\in{\rm CNF}(2,k)\\|\rho|\le2k-2}}\sum_{\substack{ \alpha = [\alpha_1,\alpha_2]^\top \in ([n]^k_{\ne})^2
      \\\sqcap (\alpha )=\rho}}|\cov(f(\mathbf{X}_{\alpha_1},\mathbf{Y}_{\alpha_1}),f(\mathbf{X}_{\alpha_2},\mathbf{Y}_{\alpha_2}))|
  \\
  \nonumber
  &\leq \sum_{\substack{\rho\in{\rm CNF}(2,k)\\|\rho|\le2k-2}}\frac{n!}{(n-|\rho|)!}2\|f\|_{\infty}^2
  \\
  \nonumber
  &\leq 2\|f\|_{\infty}^2n^{2k-2}|{\rm CNF}(2,k)|\nonumber
  \\
  \label{decp3}
  &\leq 2\|f\|_{\infty}^2n^{2k-2}2^kk!.
\end{align} 
 Finally, we bound $\Var [S_{n,k}(f)]$ from below
 using \eqref{decp1}, \eqref{decp2} and \eqref{decp3}, 
 as 
\begin{align*}
  \kappa_2(S_{n,k}(f))&= \sum_{\substack{\rho\in{\rm CNF}(2,k)\\|\rho|=2k-1}}\sum_{\substack{
      \alpha = [\alpha_1,\alpha_2]^\top \in ([n]^k_{\ne})^2\\
      \sqcap (\alpha )=\rho}}\cov (f(\mathbf{X}_{\alpha_1},\mathbf{Y}_{\alpha_1}),f(\mathbf{X}_{\alpha_2},\mathbf{Y}_{\alpha_2}))+R
  \\
  &\ge \frac{n!}{(n-2k+1)!}\Var\left[
    \sum_{\ell=1}^kf_{(\ell)}(X_1)\right]
  -2\|f\|_{\infty}^2n^{2k-2}2^kk!\\
  &= \frac{n!}{(n-2k+1)!}\left(
  \Var\left[
    \sum_{\ell=1}^kf_{(\ell)}(X_1)\right]
  -\frac{2\|f\|_{\infty}^2n^{2k-2}2^kk!}{n! / (n-2k+1)!}\right)
  \\
  &\ge \frac{n!}{(n-2k+1)!}\left(
  \Var\left[
    \sum_{\ell=1}^kf_{(\ell)}(X_1)\right]
  -\frac{2\|f\|_{\infty}^2n^{2k-2}2^kk!}{(n/2)^{2k-1}}\right), 
\end{align*}
which yields \eqref{lowb-2} by choosing
$$C(f,k):= \frac{1}{2}
\Var\left[
  \sum_{\ell=1}^kf_{(\ell)}(X_1)\right],
$$
for
$$
n
\geq 
N(f,k): = \frac{2^{3k+1}\|f\|_{\infty}^2k!}{\Var \big[
    \sum_{\ell=1}^kf_{(\ell)}(X_1)\big]
  }
.
$$
\end{Proof}
 As a consequence of Theorem~\ref{jkldd12}, 
 we have the following result. 
\begin{corollary}
  \label{cjkfl0} 
  Let $f:\mathcal{S}^k\times \mathcal{M}^{(k-1)k/2}\to\R$ be a measurable function with $k\geq 2$,
  such that $f$ satisfies Assumption~\ref{assu1}. 
  Then, for $n\ge 4(k-1)$, 
  the $j$-$th$ cumulant of the normalized $U$-statistics 
  $$
  \widebar{S}_{n,k}(f):=
        \frac{S_{n,k}(f)-\kappa_1(S_{n,k}(f))}{\sqrt{\kappa_2(S_{n,k}(f))}}
    $$ 
        satisfies
 the bound 
$$ 
    \kappa_j\left(
    \widebar{S}_{n,k}(f)
    \right)\leq
    \frac{(j!)^{k+1}}{\big(n \widetilde{C}(f,k)\big)^{j/2-1}},
$$
    $j\geq 3$, where $\widetilde{C}(f,k)>0$ depends only on
    $f$ and on $k \geq 2$.
    \end{corollary}
\begin{Proof}
 Combining the inequalities 
\begin{equation}
\label{fjklf11} 
    \frac{n!}{(n-2k+1)!}\ge(n-2k+2)^{2k-1}
     =\left(1-\frac{2k-2}n\right)^{2k-1}n^{2k-1}
     \ge
     \left(\frac{n}{2}\right)^{2k-1}, 
\end{equation}
 $n\ge4(k-1)$, and %
\begin{equation} 
 \label{djkld1111} 
  j^{j-1}<\frac{(j-1)!}{\sqrt{2\pi j}}e^j< e^jj!,
  \quad j\geq 1, 
\end{equation}
 with the bounds \eqref{fjkldf} and \eqref{lowb-2} 
for $j\ge3$, we have 
    \begin{align}
      \nonumber
      \kappa_j\left(
      \widebar{S}_{n,k}
      \right) & = \frac{\kappa_j\left(S_{n,k}(f)\right)}{\kappa_2(S_{n,k}(f))^{j/2}}
      \\
      \nonumber
      &\leq
      \Vert f\Vert_\infty^j j^{j-1} (j!)^k (k!)^{j-1}
      \frac{n^{1+(k-1)j}
        }{\left( n!C(f,k) / (n-2k+1)! \right)^{j/2}}
      \\
      \nonumber
      &\leq
      \Vert f\Vert_\infty^j j^{j-1} (j!)^k (k!)^{j-1}
      \frac{n^{1+(k-1)j}
        }{\left(C(f,k)(n/2)^{2k-1}\right)^{j/2}}
      \\
      \nonumber
      &= \Vert f\Vert_\infty^j\frac{j^{j-1}(j!)^{k}}{n^{(j-2)/2}}\frac{(k!)^{j-1}}{C(f,k)^{j/2}}
      \\
      \nonumber
      &\leq \Vert f\Vert_\infty^j
      \frac{(j!)^{k+1}}{n^{(j-2)/2}}\frac{(k!)^{j-1}e^j}{C(f,k)^{j/2}}
      \\
      \nonumber
      & \leq \frac{(j!)^{k+1}}{\big(n \widetilde{C}(f,k)\big)^{j/2-1}},
    \end{align}
    where $\widetilde{C}(f,k)>0$ depends only on
    $f$ and on $k \geq 2$.
\end{Proof}
 From Corollary~\ref{cjkfl0},
 we check that the cumulants of the normalized
 $U$-statistics
 $\widebar{S}_{n,k}$ satisfy the Statulevi\v{c}ius growth condition 
\eqref{Statuleviciuscond2}
with $\gamma := k$ and 
$$  \Delta_n := 
      \big(n\widetilde{C}(f,k)\big)^{1/2}, 
$$ 
      where $\widetilde{C}(f,k)>0$ depends only on
      $f$ and on $k \geq 2$.
 Therefore, from Proposition~\ref{l1}-$i)$
 we have the following results.
\begin{corollary}
\label{thm4.2}
(Kolmogorov bound). 
  Let $f:\mathcal{S}^k\times \mathcal{M}^{(k-1)k/2}\to\R$ be a measurable function satisfying Assumption~\ref{assu1}, with $k\geq 2$.
    For $n\ge4(k-1)$, we have 
    \begin{equation}
      \nonumber 
      \sup_{x\in\R}\left|
      \IP\left(
      \widebar{S}_{n,k}(f)     \leq x\right)-\Phi(x)\right|\leq \frac{C(f,k)}{n^{1/(2+4k)}},
    \end{equation}
    where $C(f,k)>0$ depends only on $f$ and on $k \geq 2$. 
\end{corollary}
By Corollary~\ref{cjkfl0} and 
Proposition~\ref{l1}-$ii)$,
 we also have the following result.
\begin{corollary}
\label{corg2-0}
 (Moderate deviation principle).
   Let $f:\mathcal{S}^k\times \mathcal{M}^{(k-1)k/2}\to\R$ be a measurable function satisfying Assumption~\ref{assu1}, with $k\geq 2$.
  Let $( a_n )_{n \geq 1}$ be a sequence of real numbers tending to infinity, and such that 
$ 
  a_n \ll n^{1/(2+4 k)}$. 
  Then, $\big(a_n^{-1} \widebar{S}_{n,k}(f) \big)_{n \ge1}$ satisfies a moderate deviation principle with speed $a_n^2$ and rate function $x^2/2$. 
\end{corollary} 
      \section{Upper bounds on subgraph count cumulants}
      \label{s5}
\noindent
From this section, we focus on the subgraph count in the binomial RCM
$G_H(\mathcal{X}_n)$.
 Let $\mathcal{X}_n=\{X_1,\dots,X_n\}$ be a set of i.i.d. random points on the carrier space $\mathcal{S}=\R^d$ with a common continuous distribution $\mu$. A symmetric measurable function $H:\R^d\times\R^d\to[0,1]$ is called a connection function. We consider the binomial random-connection model generated as follows: given $\mathcal{X}_n$,  we connect
 each pair of points/nodes $X_i ,X_j \in\mathcal{X}_n$, $1\leq i \ne j \leq n$, independently with probability $p_n H(X_i ,X_j )$, $0<p_n<1$. The resulting random graph is denoted by $G_H(\mathcal{X}_n)$.

 \medskip

 Given $G=(V_G,E_G)$ a connected graph with $v(G)=k\geq 2$ vertices, let 
  $$
  I_\beta :=\prod_{(\imath ,\jmath )\in E_G}\bone_{\{X_{\beta (\imath )}\sim X_{\beta (\jmath )}\}}, 
  \quad \beta \in[n]^k_{\ne},  
  $$
  where $X_\imath \sim X_\jmath $ indicates there is an edge between $X_\imath$ and $X_\jmath$,
  i.e. 
  $I_\beta =1$ 
   if and only if the
   graph with vertices $X_{\beta (1)},\ldots , X_{\beta (k)}$
   in the binomial RCM $G_H(\mathcal{X}_n)$ is isomorphic to $G$.

  \medskip

   We let
\begin{equation}
\label{defcount}
 N_G:=\frac1{a(G)}\sum_{\beta \in[n]^k_{\ne}}
 I_\beta %
\end{equation}
 denote the number of injective subgraphs of $G_H(\mathcal{X}_n)$
 isomorphic to $G$, where $X_i \sim X_j$ indicates
 that there is an edge between $X_i$ and $X_j$. 

 \medskip 

 In order to incorporate with generalized $U$-statistics in Section~\ref{s3}, we can also write the subgraph count $N_G$ as 
 \begin{equation*}
  N_G=S_{n,k}(f)=\sum_{\beta \in[n]^k_{\ne}}f\big(X_{\beta (1)},\dots,X_{\beta (k)},Y_{\beta (1),\beta (2)},\dots,Y_{\beta (1),\beta (k)},Y_{\beta (2),\beta (3)},\dots,Y_{\beta (k-1),\beta (k)}\big),
 \end{equation*}
with $f:(\R^d)^k\times [0,1]^{(k-1)k/2}\to\R$ given by 
\begin{equation}
  \label{fun1}
f(x_1,\dots,x_k,y_{1,2},\dots,y_{(k-1),k}):=\frac1{a(G)}\prod_{(\imath ,\jmath )\in E_G}\bone_{\{y_{\imath ,\jmath }\leq p_n H(x_\imath ,x_\jmath )\}},
\end{equation}
 where $\{Y_{i ,j }\}_{1\leq i < j \leq n}$ are i.i.d.
 uniform random variables on $\mathcal{M}=[0,1]$.

 \medskip
  
 Next, using
 dependency graphs and the convex analysis of planar diagrams,
 we consider the cumulant growth of the
 subgraph count $N_G$ when $p_n=o(1)$.
  Before we proceed, we need to introduce some notation.
 Since the random variables $X_1,\dots,X_n$ are i.i.d. and
 edges
 are added independently
 according to
 $$
 \bone_{\{X_i \sim X_j \}}
 =
 \bone_{\{Y_{i ,j }\leq p_n H(X_i ,X_j )\}},
\quad 1\leq i < j \leq n,
$$
 the sets of random vectors
 $$
 \{
 (I_{\alpha_1},\dots,I_{\alpha_j})
 :  \sqcap \big( [\alpha_1,\dots,\alpha_j]^\top \big) =\rho\}$$
 are identically distributed
 and have the same dependency structure
 for every $\rho\in\Pi([j] \times [k])$,
 which justifies the following definition.
\begin{definition}
 Let $j, k \geq 1$. 
 Given $\rho \in \Pi ([j]\times [k])$, we let
$$
 \kappa(\mathbf{I}_\rho):=\kappa(I_{\alpha_1},\dots,I_{\alpha_j}),
 $$
 for 
 any element
 $[\alpha_1,\ldots , \alpha_j]^\top $
 of $[n]^k_{\ne}\times \cdots \times [n]^k_{\ne}$ such that $\sqcap ( \alpha )=\rho$.
\end{definition}

\noindent
 As a consequence of \eqref{exp-1}, we have 
\begin{equation} 
\label{mean-1}
  \kappa_1(N_G)=
  \frac{n!}{(n-k)!}
  \frac{p_n^{e(G)}}{a(G)}\int_{(\R^d)^k}\prod_{(\imath ,\jmath )\in E_G}H(x_\imath ,x_\jmath )\mu ( \mathrm{d}x_1 ) \cdots\mu ( \mathrm{d}x_k ),  
\end{equation}
while the growth of higher order cumulants can be controlled as
follows. 
See Example~6.19 in \cite{JLR} and
Proposition~10.1.2 and \S~10.3.3 of \cite{feray16} for
related cumulant bounds in the
 Erd{\H o}s-R\'enyi model. 
\begin{thm}
\label{fjkdl243}
 Let $G=(V_G,E_G)$ be a connected graph with $v(G)=k\geq 2$ vertices. 
 We have the cumulant bound 
\begin{equation} 
  \label{cumbound1}
  \vert\kappa_j(N_G)\vert \leq \frac{j^{j-1}}{a(G)^j}
  \sum_{r=k}^{1+(k-1)j} |\mathcal{C}(j,k,r)|
  n^r p_n^{d(j,k,r)},
   \qquad j\geq 2, 
\end{equation}
where $\mathcal{C}(j,k,r)$ represents the collection of all partitions $\rho$ on $[j]\times[k]$ connected, non-flat and having precisely $r$ blocks, and 
\begin{equation}
  \label{min-edge}
 d(j,k,r) :=\min
 \big\{e(\rho_G) \ \!  : \ \!  \rho\in{\rm CNF}(j,k), \ |\rho|=r\big\}.
\end{equation}
\end{thm}
\begin{Proof}
  The random variables
  $$
  I_{\beta_1} :=\prod_{(\imath ,\jmath )\in E_G}\bone_{\{X_{\beta_1 (\imath )}\sim X_{\beta_1 (\jmath )}\}}
  \quad
  \mbox{and}
  \quad 
  I_{\beta_2}: =\prod_{(\imath ,\jmath )\in E_G}\bone_{\{X_{\beta_2 (\imath )}\sim X_{\beta_2 (\jmath )}\}},
  \quad \beta_1, \beta_2 \in[n]^k_{\ne}, 
  $$ 
  are independent if 
  $$
  \big\{\beta_1 (i) : i\in[k]\big\}\cap\big\{\beta_2 (i) : i\in[k]\big\}=\emptyset.$$
  Hence, by Property~C in \cite[Page~29]{MalyshevMinlos91} we have 
\begin{align} 
  \nonumber
  \kappa_j(N_G)&= \frac1{a(G)^j}\sum_{\alpha_1,\dots,\alpha_j\in[n]^k_{\ne}}\kappa(I_{\alpha_1},\dots,I_{\alpha_j})
  \\
&= \frac1{a(G)^j}\sum_{\substack{\alpha_1,\dots,\alpha_j\in[n]^k_{\ne}\\
    \mathrm{connected}}}\kappa(I_{\alpha_1},\dots,I_{\alpha_j})
\nonumber 
\\
\nonumber %
&= \frac1{a(G)^j}\sum_{ [\alpha_1,\ldots , \alpha_j]^\top \in\Lambda_j (n,k)} 
\kappa(I_{\alpha_1},\dots,I_{\alpha_j})%
, \qquad j\geq 2,
\end{align} 
where 
\begin{align*} 
\nonumber %
\Lambda_j (n,k)
&:=\left\{ \alpha = [\alpha_1,\ldots , \alpha_j]^\top \in [n]^k_{\ne}\times \cdots \times [n]^k_{\ne}
 ~\mathrm{connected}\right\}\\
&\textcolor{white}{:}
=\left\{
\alpha = [\alpha_1,\ldots , \alpha_j]^\top \in [n]^k_{\ne}\times \cdots \times [n]^k_{\ne}
 : \sqcap(\alpha) \in{\rm CNF}(j,k)\right\}. 
\end{align*}
 Now, we have 
\begin{align}
  \nonumber
  \kappa_j(N_G)&= \frac1{a(G)^j}\sum_{\rho\in{\rm CNF}(j,k)}\sum_{\substack{
 [\alpha_1,\ldots , \alpha_j]^\top \in\Lambda_j (n,k)
 \\\sqcap (
          [\alpha_1,\ldots , \alpha_j]^\top ) 
       =\rho}}\kappa(I_{\alpha_1},\dots,I_{\alpha_j}) %
    \\
  \nonumber
  &= \frac1{a(G)^j}\sum_{\rho\in{\rm CNF}(j,k)}C_n(\rho)\kappa(\mathbf{I}_\rho)
  \\
    \nonumber %
    &= \frac1{a(G)^j}\sum_{\rho\in{\rm CNF}(j,k)}
    \frac{n!}{(n-|\rho|)!}\kappa(\mathbf{I}_\rho)
  \\
  \label{cumeq-1}
  &= \frac1{a(G)^j}\sum_{r=k}^{1+(k-1)j}
  \frac{n!}{(n-r)!}
  \sum_{\substack{\rho\in{\rm CNF}(j,k)\\|\rho|=r}}\kappa(\mathbf{I}_\rho),
\end{align} 
 hence %
\begin{align}
  \label{cumu-02}
  |\kappa_j(N_G)|&\leq \frac1{a(G)^j}\sum_{r=k}^{1+(k-1)j}\frac{n!}{(n-r)!}|\mathcal{C}(j,k,r)|\max_{\rho \in\mathcal{C}(j,k,r)}|\kappa(\mathbf{I}_\rho )|. 
\end{align} 
Due to the cumulant-moment relation~\eqref{cummom1}, we have, for each $\alpha = [\alpha_1,\dots,\alpha_j]^\top \in\Lambda_{j}(n,k)$, 
\begin{align*}
  &   \kappa(I_{\alpha_1},\dots,I_{\alpha_j}) =
  \sum_{s\geq 1}
  \sum_{\sigma=\{b_1,\dots,b_s\}\in\Pi([j])}(-1)^{s-1}(s-1)!\E\left[
   \prod_{i\in b_1} I_{\alpha_i}
   \right]\times\cdots\times\E\left[
    \prod_{i\in b_s} I_{\alpha_i}
    \right]\\
  &\qquad = \E\left[\prod_{i=1}^jI_{\alpha_i}\right]+
  \sum_{s\geq 2}
  \sum_{\substack{\sigma=\{b_1,\dots,b_s\}\in\Pi([j])}}(-1)^{s-1}(s-1)!\E\left[
    \prod_{i\in b_1} I_{\alpha_i}
    \right]\times\cdots\times\E\left[
    \prod_{i\in b_s} I_{\alpha_i}
    \right].
\end{align*}
Moreover, denoting by
 $\sqcap (\alpha )_G$ the connected graph built on
$\sqcap (\alpha )\in\Pi([j] \times [k])$
in Definition~\ref{def:diagram-graph}, we have
\begin{eqnarray*}
  \E\left[\prod_{i=1}^jI_{\alpha_i}\right]&=&
  \E\left[\prod_{i=1}^j\left(\prod_{(\imath,\jmath)\in E(G)}\bone_{\{X_{\alpha_i(\imath)}\sim X_{\alpha_i(\jmath)}\}}\right)\right]\\
  &=& p_n^{e(\sqcap (\alpha  )_G)}\int_{(\R^d)^{|\sqcap (\alpha  )|}}\prod_{(\imath ,\ell)\in E(\sqcap (\alpha  )_G)}H(x_\imath ,x_\ell)\mu(\mathrm{d}x_1)\cdots\mu(\mathrm{d}x_{|\sqcap (\alpha )|}) \\
  &\leq &p_n^{e(\sqcap (\alpha )_G)},
\end{eqnarray*}
 where the last inequality is due to the facts that $0\leq H\le1$ and $\mu$ is a probability measure on $\R^d$.  
For each $\rho\in\Pi([j]\times[k])$ and $b\subset[j]$, we denote 
$\rho_{|b}:=\rho_{|\{\pi_i\ \!  : \ \!  i\in b\}}$ the partition of $\cup_{i\in b}\pi_i$ obtained restricting $\rho$ on $\cup_{i\in b}\pi_i$, i.e. 
$$\rho_{|b} =\{a\cap (\cup_{i\in b}\pi_i) \ \!  : \ \!  a\in\rho\}.$$
 Therefore, for any $\sigma=\{b_1,\dots,b_r\}\in\Pi([j])$ with $r\geq 2$ we have,
 letting $\rho:=\sqcap (a)$, 
\begin{align}
  \prod_{i=1}^r
  \E\left[
    \prod_{j\in b_i} I_{\alpha_j}
 \right]
  &=\prod_{i=1}^r\E\left[\prod_{\ell\in b_i}\left(\prod_{(\imath,\jmath)\in E(G)}\bone_{\{X_{\alpha_\ell(\imath)}\sim X_{\alpha_\ell(\jmath)}\}}\right)\right]\nonumber\\
  &=\prod_{i=1}^r p_n^{e((\rho_{|b_i})_G)}\int_{(\R^d)^{\left|\rho_{|b_i}\right|}}\prod_{(\imath ,\ell)\in E((\rho_{|b_i})_G)}H(x_\imath ,x_\ell)\mu(\mathrm{d}x_1)\cdots\mu(\mathrm{d}x_{| \rho_{|b_i} |})
  \nonumber\\
  &\leq p_n^{\sum_{i=1}^re((\rho_{|b_i})_G)}
  \nonumber
  \\
  & \leq p_n^{e(\rho_G)}
   \nonumber
  \\
   & = p_n^{e(\sqcap (\alpha )_G)}.  
\nonumber
\end{align} 
 Next, 
 we have
\begin{eqnarray*}
  \max_{\rho \in\mathcal{C}(j,k,r)}|\kappa(\mathbf{I}_\rho )|&\leq &p_n^{d(j,k,r)}
  \sum_{s=1}^j
  \sum_{\sigma=\{b_1,\dots,b_s\}\in\Pi([j])}(s-1)!\nonumber\\
&=&p_n^{d(j,k,r)}\sum_{s=1}^j(s-1)!S(j,s)\nonumber\\
&\leq & \frac{p_n^{d(j,k,r)}}{j}
\sum_{s=1}^j(j)_sS(j,s)\nonumber\\
&\leq &p_n^{d(j,k,r)}j^{j-1}, 
\end{eqnarray*}
where $S(j,r)$ is the Stirling number of the second kind.
Therefore, together with \eqref{cumu-02}, we obtain 
\begin{eqnarray*}
|\kappa_j(N_G)|&\leq &
\frac{1}{a(G)^j} \sum_{r=k}^{1+(k-1)j}
\frac{n!}{(n-r)!}
|\mathcal{C}(j,k,r)| p_n^{d(j,k,r)}j^{j-1} \\
&\leq &\frac{j^{j-1}}{a(G)^j}\sum_{r=k}^{1+(k-1)j} |\mathcal{C}(j,k,r)|n^r p_n^{d(j,k,r)},
\end{eqnarray*}
 which proves \eqref{cumbound1}.
\end{Proof}
\section{Variance lower bounds for subgraph counts} 
\label{s6}
 \noindent
To prove Theorem~\ref{th7.5}, we shall derive a variance lower bound for the subgraph count $N_G$ in the binomial RCM. For this purpose, we need to introduce some notation. Recall $f$ from \eqref{fun1}, we further denote, for $1\leq \ell\leq k$,
\begin{align}
  \nonumber
  f_{(\ell)}(x)&:=\frac1{a(G)}\E\left[\prod_{(\imath ,\jmath )\in E_G}\bone_{\left\{Y_{\imath ,\jmath }\leq p_n H(X_\imath ,X_\jmath )\right\}}\middle|X_{\ell}=x\right]
  \\
  \nonumber
  &\textcolor{white}{:}
   =\int_{(\R^d)^{k-1}\times [0,1]^{(k-1)k/2}}f(\mathbf{x},\mathbf{y})\mu^{\otimes  (k-1)}
  ( \mathrm{d}x_1 , \ldots , \mathrm{d}x_{\ell-1} , \mathrm{d}x_{\ell+1} ,
   \ldots , \mathrm{d}x_k ) \mathrm{d}\mathbf{y}
   \\
  \label{fun2}
  &\textcolor{white}{:}
  =\frac{p_n^{e(G)}}{a(G)}\int_{(\R^d)^{k-1}}\prod_{(\imath ,\jmath )\in E(G)}H(x_\imath ,x_\jmath )\mu^{ \otimes ( k-1 )}
  ( \mathrm{d}x_1, \ldots , \mathrm{d}x_{\ell-1} , \mathrm{d}x_{\ell+1},
   \ldots , \mathrm{d}x_k ),
\end{align}
where $\mathbf{x}=(x_1,\dots,x_{\ell-1},x,x_{\ell+1},\dots,x_k )$.

\begin{prop}
  \label{p6.1}
  Suppose that the
  function $f$ in \eqref{fun1}
  satisfies Assumption~\ref{assu1}.
   Then, we have the following variance lower bound for subgraph counts: 
   \begin{equation}
     \label{var-low}
    \Var [ N_G ]\ge\frac{C}{a(G)^2}
    \left(
    \frac{n!}{(n-2k+1)!}p_n^{2e(G)}+\frac{n!}{(n-k)!}p_n^{e(G)}\right), 
  \end{equation}
  where $C>0$ is a constant independent of $n\geq 2k-1$.
\end{prop}
\begin{Proof}
 We reuse the notation 
$$
 I_\beta :=\prod_{(\imath ,\jmath )\in E_G}\bone_{\{X_{\beta (\imath )}\sim X_{\beta (\jmath )}\}},
 \quad \beta \in[n]^k_{\ne}.
 $$ 
Because the indicator random variable $I_\beta $, $\beta \in[n]^k_{\ne}$, indicates that there is a copy of $G$ on $\big\{X_{\beta (1)},\dots,X_{\beta (k)}\big\}$, 
it is easy to see that for any $\alpha_1,\alpha_2\in[n]^k_{\ne}$ we have  
\begin{eqnarray*}
  \E[I_{\alpha_1} \mid I_{\alpha_2}=1]\ge\E[I_{\alpha_1}],
\end{eqnarray*}
 hence 
$$
 \cov(I_{\alpha_1},I_{\alpha_2})\geq 0, \quad
  \alpha_1,\alpha_2\in[n]^k_{\ne}. 
$$ 
  This further implies $\kappa(\mathbf{I}_\rho)\geq 0$
  for any $\rho\in{\rm CNF}(2,k)$.
  For any $\rho\in{\rm CNF}(2,k)$ with $|\rho|=2k-1$, there is only one block containing exactly two elements in $\rho$ and the rest $2k-2$ blocks containing exactly one element. Without loss of generality, we assume $\rho\in{\rm CNF}(2,k)$ with $|\rho|=2k-1$ and $\{(1,1),(2,t)\}\in\rho$
  for some $t\in [k]$. Then, we have
\begin{eqnarray*}
    \kappa(\mathbf{I}_\rho)&=&\E\left[\left(\prod_{(\imath ,\jmath )\in E_G}\bone_{\{X_{\imath }\sim X_{\jmath }\}}\right)\left(\prod_{(\imath ,\jmath )\in E_G}\bone_{\{X'_\imath \sim X'_\jmath \}}\right)\right]-\left(\E\left[\prod_{(\imath ,\jmath )\in E_G}\bone_{\{X_\imath \sim X_\jmath \}}\right]\right)^2\nonumber\\
    &=&a(G)^2\big(\E[f_{(1)}(X_1)f_{(t)}(X_1)]-\E[f_{(1)}(X_1)]^2\big),
  \end{eqnarray*}
  where $\{X_1',\dots,X_k'\}$ are i.i.d. copies of $X\sim\mu$, independent of $\{X_1,\dots,X_k\}$, with $X_t'$ being replaced by $X_1$. From Assumption~\ref{assu1}, we have 
  \begin{eqnarray*}
    \sum_{\substack{\rho\in{\rm CNF}(2,k)\\|\rho|=2k-1}}\kappa(\mathbf{I}_\rho)&=&a(G)^2\left(\E\left[\left(\sum_{i=1}^kf_{(i)}(X_1)\right)^2\right]-k^2\E[f_{(1)}(X_1)]^2\right)\nonumber\\
    &=&a(G)^2\Var\left[\sum_{i=1}^kf_{(i)}(X_1)\right]>0.
  \end{eqnarray*}
  On the other hand, combining the above with \eqref{fun2}, we find that 
     \begin{align}
   \nonumber
      & \sum_{\substack{\rho\in{\rm CNF}(2,k)\\|\rho|=2k-1}}\kappa(\mathbf{I}_\rho)
      =
      a(G)^2\Var\left[\sum_{i=1}^kf_{(i)}(X_1)\right]
      \\
   \nonumber
            & \quad =
      p_n^{2e(G)}
   \Var\left[\sum_{i=1}^kf_{(i)}
\int_{(\R^d)^{k-1}}\prod_{(\imath ,\jmath )\in E(G)}H(x_\imath ,x_\jmath )\mu^{ \otimes ( k-1 )}
  ( \mathrm{d}x_1, \ldots , \mathrm{d}x_{i-1} , \mathrm{d}x_{i+1},
\ldots , \mathrm{d}x_k )
\right]
   \\
\label{varlbound1}
         & \quad     =C_1 p_n^{2e(G)},
    \end{align}
    where
    \begin{align*} 
     C_1: & = 
     C_1(\mu,H)
     \\
      & = \Var\left[\sum_{i=1}^kf_{(i)}
\int_{(\R^d)^{k-1}}\prod_{(\imath ,\jmath )\in E(G)}H(x_\imath ,x_\jmath )\mu^{ \otimes ( k-1 )}
  ( \mathrm{d}x_1, \ldots , \mathrm{d}x_{i-1} , \mathrm{d}x_{i+1},
\ldots , \mathrm{d}x_k )
\right]>0
   \end{align*} 
 is a constant independent of $n \geq 1$ 
   and 
   $\mathbf{x}=(x_1,\dots,x_{\ell-1},X_1,x_{\ell+1},\dots,x_k )$.

  \medskip
   
  Taking $j=2$ in \eqref{cumeq-1}, together with \eqref{varlbound1}, we obtain
  the following lower bound for the variance of subgraph counts: 
\begin{align*}
  \nonumber %
  \Var [ N_G ]&=\frac1{a(G)^2}\sum_{r=k}^{2k-1}
  \frac{n!}{(n-r)!}
  \sum_{\substack{\rho\in{\rm CNF}(2,k)\\|\rho|=r}}\kappa(\mathbf{I}_\rho)\nonumber\\
  &\ge\frac1{a(G)^2}\left(
  \frac{n!}{(n-k)!}
  \sum_{\substack{\rho\in{\rm CNF}(2,k)\\|\rho|=k}}\kappa(\mathbf{I}_\rho)+
  \frac{n!}{(n-2k+1)!}
  \sum_{\substack{\rho\in{\rm CNF}(2,k)\\|\rho|=2k-1}}\kappa(\mathbf{I}_\rho)\right)
  \nonumber\\
  &\ge\frac1{a(G)^2}\left(
  \frac{n!}{(n-k)!}
  \kappa(\mathbf{I}_{\tilde{\rho}}) + %
  \frac{n!}{(n-2k+1)!} 
  \sum_{\substack{\rho\in{\rm CNF}(2,k)\\|\rho|=2k-1}}\kappa(\mathbf{I}_\rho)\right)
  \nonumber,
\end{align*} 
 where $\tilde{\rho}$ is the partition
 $\tilde{\rho} =
 \big\{
 \{(1,1),(2,1)\},
 \ldots  ,
 \{(1,k),(2,k)\}
 \big\}$
 of $[2]\times [k]$
 Next, for any $\beta = [ \beta_1, \beta_2]^\top 
 \in ([n]^k_{\ne})^2 %
 $ such that $\sqcap (\beta ) = \tilde{\rho}$,
 we have 
  $$
 \kappa(\mathbf{I}_{\tilde{\rho}}) 
 =
 \Var [ I_\beta ]
  =
  \theta p_n^r (1-\theta p_n^r )
  =
  \theta p_n^{e(G)}\left(1-\theta p_n^{e(G)}\right)
  $$
  where
  $$\theta:=\int_{(\R^d)^k}\prod_{(\imath ,\jmath )\in E(G)}H(x_\imath ,x_\jmath )
\mu^{\otimes k} ( \mathrm{d}x_1 , \ldots , \mathrm{d}x_k),
$$
    hence 
\begin{align*} 
 \Var [ N_G ] & \geq
    \frac1{a(G)^2}\left(
  \frac{n!}{(n-k)!}
  \theta p_n^{e(G)}\left(1-\theta p_n^{e(G)}\right)+
  \frac{n!}{(n-2k+1)!}C_1 p_n^{2e(G)}\right)
  \nonumber\\
  &\ge\frac{C_2}{a(G)^2}\left(
  \frac{n!}{(n-2k+1)!}p_n^{2e(G)}+\frac{n!}{(n-k)!}p_n^{e(G)}\right), 
\end{align*} 
  with $C_2:=\min \big( C_1 ,
\theta \big(
1-\theta p_n^{e(G)}
\big)\big)$.
\end{Proof}
\section{Growth rates of subgraph count cumulants}
\label{s7}
\noindent
The general cumulant upper bound of the subgraph count $N_G$ established in
Theorem~\ref{fjkdl243} does not yield an 
explicit asymptotic growth order. 
In this section, we perform a more detailed analysis
of cumulant growth rates by identifying
the leading terms of the form $n^rp_n^{d(j,k,r)}$ in \eqref{cumbound1},
$r=k,\dots,1+(k-1)j$, where $d(j,k,r)$ is defined in \eqref{min-edge},
 with in particular 
\begin{equation}
  \label{end1}
  d_k=e(G), \quad d_{1+(k-1)j}=je(G).
\end{equation}
To this end, we start by deriving cumulant growth rates via the convex analysis of planar diagrams used in \cite{LiuPrivault24}
 in the random-connection model and introduced in 
 \cite{JLR} for the Erd{\H o}s-R\'enyi model. 
 Firstly, we adopt some notation from \cite{LiuPrivault24}. 
\begin{definition}
Let $G$ be a connected graph with $k\ge2$ vertices. For $j\ge1$, we let 
\begin{enumerate}[i)]
\item
  $\displaystyle \Sigma_j(G):=\big\{(x(\rho_G),y(\rho_G)):=(jk-v(\rho_G),je(G)-e(\rho_G)): \rho\in{\rm CNF}(j,k)\big\}$,
   \\ 
where for each $\rho\in\Pi([j]\times[k])$, $\rho_G$ is the graph associated to $\rho$ by Definition~\ref{def:diagram-graph}; 
\item and we let $\widehat{\Sigma}_j(G)$ denote the upper boundary of the convex hull of $\Sigma_j(G)$. 
\end{enumerate}
\end{definition}
\begin{remark}\label{end2}
  According to \eqref{end1}, we can see that two endpoints of
  the upper boundary of the convex hull of $\Sigma_j(G)$
  have coordinates $(j-1,0)$ and $((j-1)k,(j-1)e(G))$.
\end{remark}
 The following example was considered in \cite[Example~5.2]{LiuPrivault24} in the
 Poisson random-connection model.
\begin{example}
  Let $G=C_3$ be a triangle, i.e. $v(G)=e(G)=3$.
  We have 
  $$
 \left\{
 \begin{array}{l}
  \Sigma_2(C_3) = \{(3,3),(2,1),(1,0)\},
   \medskip
   \\ 
  \Sigma_3(C_3) = \{(6, 6), (5, 4), (4, 3), (5, 3), (4, 2), (4, 1), (3, 1), (3, 0), (2, 0)\},
     \medskip
   \\ 
   \Sigma_4(C_3) = \{(9, 9), (8, 7), (7, 6), (8, 6), (7, 5), (7, 4), (6, 4), (6, 3),(5, 3), (7, 3),
   \smallskip
   \\
   \qquad \qquad \ \ \ \!  
   (6, 2), (5, 2), (7, 2), (6, 1), (5, 1), (4, 1), (6, 0), (5, 0), (4, 0), (3, 0)\},
 \end{array}
 \right.
   $$ 
 see Figure~\ref{fig1}-$(a)$ 
 with $j=3$, $k=3$, $e(G)=3$, and
 Figure~\ref{fig1}-$(b)$
 $j=4$, $k=3$, $e(G)=3$.

\begin{figure}[H]
  \begin{subfigure}{.5\textwidth}
    \vskip1cm
  \centering
\begin{tikzpicture}[scale=0.9] 
\draw[step={(1cm,0.5cm)}, very thin, gray!40] (0,0) grid (7,3.5);
\draw[->] (0,0) -- (0,3.5);
\node[anchor=east,font=\small] at (-0.05,3.5) {$y$};
\draw[->] (0,0) -- (7,0);
\node[anchor=north,font=\small] at (7,-0.05) {$x$};
\draw[thick] (-0.05,0.5) -- (0.05,0.5);
\node[anchor=east,font=\small] at (-0.05,0.5) {1};
\draw[thick] (-0.05,1) -- (0.05,1);
\node[anchor=east,font=\small] at (-0.05,1) {2};
\draw[thick] (-0.05,2) -- (0.05,2);
\node[anchor=east,font=\small] at (-0.05,2) {4};
\draw[thick] (-0.05,3) -- (0.05,3);
\node[anchor=east,font=\small] at (-0.05,3) {6};
\draw[thick] (-0.05,1.5) -- (0.05,1.5);
\node[anchor=east,font=\small] at (-0.05,1.5) {3};
\draw[thick] (-0.05,2.5) -- (0.05,2.5);
\node[anchor=east,font=\small] at (-0.05,2.5) {5};
\node[anchor=east,font=\small] at (-0.05,3) {6};
\draw[thick] (1,-0.05) -- (1,0.05);
\node[anchor=north,font=\small] at (1,-0.05) {1};
\draw[thick] (2,-0.05) -- (2,0.05);
\node[anchor=north,font=\small] at (2,-0.05) {2};
\draw[thick] (3,-0.05) -- (3,0.05);
\node[anchor=north,font=\small] at (3,-0.05) {3};
\draw[thick] (4,-0.05) -- (4,0.05);
\node[anchor=north,font=\small] at (4,-0.05) {4};
\draw[thick] (5,-0.05) -- (5,0.05);
\node[anchor=north,font=\small] at (5,-0.05) {5};
\draw[thick] (6,-0.05) -- (6,0.05);
\node[anchor=north,font=\small] at (6,-0.05) {6};
\draw[very thick,red] (2,0) -- (6,3);
\filldraw [blue] (2,0) circle (2.0pt);
\filldraw [blue] (3,0) circle (2.0pt);
\filldraw [blue] (3,0.5) circle (2.0pt);
\filldraw [blue] (4,0.5) circle (2.0pt);
\filldraw [blue] (4,1) circle (2.0pt);
\filldraw [blue] (5,1.5) circle (2.0pt);
\filldraw [blue] (4,1.5) circle (2.0pt);
\filldraw [blue] (5,2) circle (2.0pt);
\filldraw [blue] (6,3) circle (2.0pt);
    \end{tikzpicture}
\caption{Plot of $\Sigma_3(C_3)$.} %
\end{subfigure}
\begin{subfigure}{.5\textwidth}
  \centering
  \begin{tikzpicture}[scale=0.9] 
\draw[step=0.5cm, very thin, gray!40] (0,0) grid (5,5);
\draw[->] (0,0) -- (0,5);
\node[anchor=east,font=\small] at (-0.05,5) {$y$};
\draw[->] (0,0) -- (5,0);
\node[anchor=north,font=\small] at (5,-0.05) {$x$};
\draw[thick] (-0.05,0.5) -- (0.05,0.5);
\node[anchor=east,font=\small] at (-0.05,0.5) {1};
\draw[thick] (-0.05,1) -- (0.05,1);
\node[anchor=east,font=\small] at (-0.05,1) {2};
\draw[thick] (-0.05,1.5) -- (0.05,1.5);
\node[anchor=east,font=\small] at (-0.05,1.5) {3};
\draw[thick] (-0.05,2) -- (0.05,2);
\node[anchor=east,font=\small] at (-0.05,2) {4};
\draw[thick] (-0.05,2.5) -- (0.05,2.5);
\node[anchor=east,font=\small] at (-0.05,2.5) {5};
\draw[thick] (-0.05,3) -- (0.05,3);
\node[anchor=east,font=\small] at (-0.05,3) {6};
\draw[thick] (-0.05,3.5) -- (0.05,3.5);
\node[anchor=east,font=\small] at (-0.05,3.5) {7};
\draw[thick] (-0.05,4) -- (0.05,4);
\node[anchor=east,font=\small] at (-0.05,4) {8};
\draw[thick] (-0.05,4.5) -- (0.05,4.5);
\node[anchor=east,font=\small] at (-0.05,4.5) {9};
\draw[thick] (0.5,-0.05) -- (0.5,0.05); %
\node[anchor=north,font=\small] at (0.5,-0.05) {1}; %
\draw[thick] (1,-0.05) -- (1,0.05);
\node[anchor=north,font=\small] at (1,-0.05) {2};
\draw[thick] (1.5,-0.05) -- (1.5,0.05);
\node[anchor=north,font=\small] at (1.5,-0.05) {3};
\draw[thick] (2,-0.05) -- (2,0.05);
\node[anchor=north,font=\small] at (2,-0.05) {4};
\draw[thick] (2.5,-0.05) -- (2.5,0.05);
\node[anchor=north,font=\small] at (2.5,-0.05) {5};
\draw[thick] (3,-0.05) -- (3,0.05);
\node[anchor=north,font=\small] at (3,-0.05) {6};
\draw[thick] (3.5,-0.05) -- (3.5,0.05);
\node[anchor=north,font=\small] at (3.5,-0.05) {7};
\draw[thick] (4,-0.05) -- (4,0.05);
\node[anchor=north,font=\small] at (4,-0.05) {8};
\draw[thick] (4.5,-0.05) -- (4.5,0.05);
\node[anchor=north,font=\small] at (4.5,-0.05) {9};
\draw[very thick,red] (1.5,0) -- (4.5,4.5);
\filldraw [blue] (1.5,0) circle (2.0pt);
\filldraw [blue] (2,0) circle (2.0pt);
\filldraw [blue] (2.5,0) circle (2.0pt);
\filldraw [blue] (3,0) circle (2.0pt);
\filldraw [blue] (2,0.5) circle (2.0pt);
\filldraw [blue] (2.5,0.5) circle (2.0pt);
\filldraw [blue] (3,0.5) circle (2.0pt);
\filldraw [blue] (3.5,1) circle (2.0pt);
\filldraw [blue] (3.0,1) circle (2.0pt);
\filldraw [blue] (2.5,1) circle (2.0pt);
\filldraw [blue] (3.5,1.5) circle (2.0pt);
\filldraw [blue] (2.5,1.5) circle (2.0pt);
\filldraw [blue] (3,1.5) circle (2.0pt);
\filldraw [blue] (3,2) circle (2.0pt);
\filldraw [blue] (3.5,2) circle (2.0pt);
\filldraw [blue] (3.5,2.5) circle (2.0pt);
\filldraw [blue] (4,3) circle (2.0pt);
\filldraw [blue] (3.5,3) circle (2.0pt);
\filldraw [blue] (4,3.5) circle (2.0pt);
\filldraw [blue] (4.5,4.5) circle (2.0pt);
\end{tikzpicture}%
  \caption{Plot of $\Sigma_4(C_3)$.} %
\end{subfigure}
\caption{Set $\Sigma_n(C_3)$ and upper boundary of its convex hull (in red) for $n=3,4$.}
\label{fig1}
\end{figure}

\vspace{-0.4cm}

\noindent
From Figure~\ref{fig1} we can also check that 
$$
\begin{array}{ll}
  \widehat{\Sigma}_3(C_3) \cap
  {\Sigma}_3(C_3)
  = & \{(6, 6), (4, 3), (2, 0)\},
     \medskip
   \\ 
   \widehat{\Sigma}_4(C_3)
   \cap
   {\Sigma}_4(C_3)
   = & \{(9, 9), (7, 6), (5, 3), (3, 0)\}. 
 \end{array}
$$
\end{example}
 We recall the following definition. 
\begin{definition}{\cite{luczak}, \cite[pages~64-65]{JLR}} 
 A graph $G$ is strongly balanced if 
\begin{equation}
     \label{strongbaldef}
    \frac{e(H)}{v(H)-1}\leq \frac{e(G)}{v(G)-1},
    \qquad
    H\subset G. %
  \end{equation}
  \end{definition}
\begin{thm}
  \label{th7.5}
  Let $G$ be a strongly balanced connected graph. 
\begin{enumerate}[a)]
\item If $n^{-(v(G)-1)/e(G)}\ll p_n$,
  then for $n$ large enough 
  we have 
  \begin{equation}
    \label{upp-1}
    |\kappa_j(N_G)|\leq
    j^{j-1}j!^{v(G)} \frac{v(G)!^{j-1}}{a(G)^j}n^{1+(v(G)-1)j}p_n^{je(G)},
    \quad j \geq 2. 
  \end{equation}
\item If $p_n\ll n^{-(v(G)-1)/e(G)}$,
  then for $n$ large enough we have 
    \begin{equation}
      \label{upp-2}
      |\kappa_j(N_G)|\leq
      j^{j-1}j!^{v(G)}
      \frac{v(G)!^{j-1}}{a(G)^j}n^{v(G)}p_n^{e(G)},
    \quad j \geq 2. 
  \end{equation}
\end{enumerate}
\end{thm}
\begin{Proof}
  We let $k:=v(G)$.
  According to \cite[Proposition~6.4]{LiuPrivault24}, if $G$ is a
   strongly balanced graph, then the upper boundary of the convex hull of $\Sigma_j(G)$, $j\ge2$,
   is a line segment, whose endpoints are $(j-1,0)$ and $((j-1)k,(j-1)e(G))$,
   see Remark~\ref{end2} and Figure~\ref{fig1-3} for $j=k=3$.
  Then, for any $r \in \{k+1,\ldots,(k-1)j\}$, the asymptotic order
  of summands appearing in the upper bound of \eqref{cumbound1} is corresponding to a point $(jk-r, je(G)-d(j,k,r))$ denoted by the blue dots in Figure~\ref{fig1-3}. 
  For every such point $(jk-r, je(G)-d(j,k,r))$ located within the convex hull, comparing the ratio of a line segment that connects itself with the right endpoint $((j-1)k,(j-1)e(G))$ of the upper boundary of the convex hull with the ratio of the upper boundary, we have 
 \begin{equation}
   \nonumber 
   \frac{(j-1)e(G)-(je(G)-d(j,k,r))}{(j-1)k-(jk-r)} \ge\frac{(j-1)e(G)}{
     (j-1)k-(j-1)} = \frac{e(G)}{k-1}, 
\end{equation}
 hence, choosing $j \geq 2$, 
 \begin{equation}
   \label{convex3}
 \frac{d(j,k,r)-e(G)}{r-k}\ge\frac{e(G)}{k-1}, 
\end{equation}
 see for example the dashed line in Figure~\ref{fig1-3}-$(a)$
   with
   $j=k=e(G)=3$, $r=5$ and $d(3,3,5) = 8$,
   in the framework of Example~\ref{fjkld243}.
 
 \vskip-1cm 
  
\begin{figure}[H]
  \begin{subfigure}{.5\textwidth}
    \vskip1cm
  \centering
\begin{tikzpicture}[scale=0.9] 
\draw[step={(1cm,0.5cm)}, very thin, gray!40] (0,0) grid (7,3.5);
\draw[->] (0,0) -- (0,3.5);
\node[anchor=east,font=\small] at (-0.05,3.5) {$y$};
\draw[->] (0,0) -- (7,0);
\node[anchor=north,font=\small] at (7,-0.05) {$x$};
\draw[thick] (-0.05,0.5) -- (0.05,0.5);
\node[anchor=east,font=\small] at (-0.05,0.5) {1};
\draw[thick] (-0.05,1) -- (0.05,1);
\node[anchor=east,font=\small] at (-0.05,1) {2};
\draw[thick] (-0.05,2) -- (0.05,2);
\node[anchor=east,font=\small] at (-0.05,2) {4};
\draw[thick] (-0.05,3) -- (0.05,3);
\node[anchor=east,font=\small] at (-0.05,3) {6};
\draw[thick] (-0.05,1.5) -- (0.05,1.5);
\node[anchor=east,font=\small] at (-0.05,1.5) {3};
\draw[thick] (-0.05,2.5) -- (0.05,2.5);
\node[anchor=east,font=\small] at (-0.05,2.5) {5};
\node[anchor=east,font=\small] at (-0.05,3) {6};
\draw[thick] (1,-0.05) -- (1,0.05);
\node[anchor=north,font=\small] at (1,-0.05) {1};
\draw[thick] (2,-0.05) -- (2,0.05);
\node[anchor=north,font=\small] at (2,-0.05) {2};
\draw[thick] (3,-0.05) -- (3,0.05);
\node[anchor=north,font=\small] at (3,-0.05) {3};
\draw[thick] (4,-0.05) -- (4,0.05);
\node[anchor=north,font=\small] at (4,-0.05) {4};
\draw[thick] (5,-0.05) -- (5,0.05);
\node[anchor=north,font=\small] at (5,-0.05) {5};
\draw[thick] (6,-0.05) -- (6,0.05);
\node[anchor=north,font=\small] at (6,-0.05) {6};
\draw[very thick,red] (2,0) -- (6,3);
\draw[very thick, purple, dashed] (4,0.5) -- (6,3);
\filldraw [blue] (2,0) circle (2.0pt);
\filldraw [blue] (3,0) circle (2.0pt);
\filldraw [blue] (3,0.5) circle (2.0pt);
\filldraw [blue] (4,0.5) circle (2.0pt);
\filldraw [blue] (4,1) circle (2.0pt);
\filldraw [blue] (5,1.5) circle (2.0pt);
\filldraw [blue] (4,1.5) circle (2.0pt);
\filldraw [blue] (5,2) circle (2.0pt);
\filldraw [blue] (6,3) circle (2.0pt);
    \end{tikzpicture}
\caption{Plot of $\Sigma_3(C_3)$.} %
\end{subfigure}
  \begin{subfigure}{.5\textwidth}
    \vskip1cm
  \centering
\begin{tikzpicture}[scale=0.9] 
\draw[step={(1cm,0.5cm)}, very thin, gray!40] (0,0) grid (7,3.5);
\draw[->] (0,0) -- (0,3.5);
\node[anchor=east,font=\small] at (-0.05,3.5) {$y$};
\draw[->] (0,0) -- (7,0);
\node[anchor=north,font=\small] at (7,-0.05) {$x$};
\draw[thick] (-0.05,0.5) -- (0.05,0.5);
\node[anchor=east,font=\small] at (-0.05,0.5) {1};
\draw[thick] (-0.05,1) -- (0.05,1);
\node[anchor=east,font=\small] at (-0.05,1) {2};
\draw[thick] (-0.05,2) -- (0.05,2);
\node[anchor=east,font=\small] at (-0.05,2) {4};
\draw[thick] (-0.05,3) -- (0.05,3);
\node[anchor=east,font=\small] at (-0.05,3) {6};
\draw[thick] (-0.05,1.5) -- (0.05,1.5);
\node[anchor=east,font=\small] at (-0.05,1.5) {3};
\draw[thick] (-0.05,2.5) -- (0.05,2.5);
\node[anchor=east,font=\small] at (-0.05,2.5) {5};
\node[anchor=east,font=\small] at (-0.05,3) {6};
\draw[thick] (1,-0.05) -- (1,0.05);
\node[anchor=north,font=\small] at (1,-0.05) {1};
\draw[thick] (2,-0.05) -- (2,0.05);
\node[anchor=north,font=\small] at (2,-0.05) {2};
\draw[thick] (3,-0.05) -- (3,0.05);
\node[anchor=north,font=\small] at (3,-0.05) {3};
\draw[thick] (4,-0.05) -- (4,0.05);
\node[anchor=north,font=\small] at (4,-0.05) {4};
\draw[thick] (5,-0.05) -- (5,0.05);
\node[anchor=north,font=\small] at (5,-0.05) {5};
\draw[thick] (6,-0.05) -- (6,0.05);
\node[anchor=north,font=\small] at (6,-0.05) {6};
\draw[very thick,red] (2,0) -- (6,3);
\draw[very thick, purple, dashed] (2,0) -- (4,0.5);
\filldraw [blue] (2,0) circle (2.0pt);
\filldraw [blue] (3,0) circle (2.0pt);
\filldraw [blue] (3,0.5) circle (2.0pt);
\filldraw [blue] (4,0.5) circle (2.0pt);
\filldraw [blue] (4,1) circle (2.0pt);
\filldraw [blue] (5,1.5) circle (2.0pt);
\filldraw [blue] (4,1.5) circle (2.0pt);
\filldraw [blue] (5,2) circle (2.0pt);
\filldraw [blue] (6,3) circle (2.0pt);
    \end{tikzpicture}
\caption{Plot of $\Sigma_3(C_3)$.} %
\end{subfigure}
\caption{Set $\Sigma_3(C_3)$ and upper boundary of its convex hull (in red).}
\label{fig1-3}
\end{figure}

\vspace{-0.4cm}

\noindent
Similarly,  
for any $r \in \{k+1,\ldots,(k-1)j\}$, when comparing ratio of line segment between point $(jk-r, je(G)-d(j,k,r))$ and
 the left endpoint of the upper boundary
 $(j-1,0)$ with the ratio of the upper boundary $\widehat{\Sigma}_j(G)$,
 if $(jk-r,je(G)-d(j,k,r))\in \widehat{\Sigma}_j(G)$, then we have 
 \begin{equation}
   \nonumber 
  \frac{je(G)-d(j,k,r)}{jk-r-(j-1)}=\frac{(j-1)e(G)}{(j-1)k-(j-1)}, 
\end{equation}
 hence 
 \begin{equation}
   \label{convex1}
  \frac{je(G)-d(j,k,r)}{jk-r-(j-1)}=\frac{e(G)}{k-1}.
\end{equation}
On the other hand, if $(jk-r,je(G)-d(j,k,r))\notin \widehat{\Sigma}_j(G)$, we have 
\begin{equation}
  \label{convex2}
  \frac{je(G)-d(j,k,r)}{jk-r-(j-1)} %
  <\frac{e(G)}{k-1},
\end{equation}
 see for example the dashed line in Figure~\ref{fig1-3}-$(b)$
   with
   $j=k=e(G)=3$, $r=5$ and $d(3,3,5) = 8$, 
   as in Example~\ref{fjkld243}. 
\begin{enumerate}[a)]
\item When $n^{-(k-1)/e(G)}\ll p_n$, 
  using \eqref{convex1} and \eqref{convex2}, 
  for any $r\in \{k,\ldots,(k-1)j\}$ we have 
  \begin{align}
    \nonumber 
    \frac{n^{1+(k-1)j}p_n^{je(G)}}{n^rp_n^{d(j,k,r)}}&= n^{1+(k-1)j-r}p_n^{je(G)-d(j,k,r)}
    \\
    \nonumber 
  &= \left( np_n^{\frac{je(G)-d(j,k,r)}{1+(k-1)j-r}}\right)^{1+(k-1)j-r}
  \\
    \nonumber 
  &\gtrsim \left( np_n^{\frac{e(G)}{k-1}}\right)^{1+(k-1)j-r}
  \\
  \label{fjkldfaa1}
  & \gg  1.
\end{align}
  Therefore, after combining with \eqref{cumbound1},
  for $n$ large enough we have 
\begin{align*}
  \vert\kappa_j(N_G)\vert &\leq  \frac{j^{j-1}}{a(G)^j}
  \sum_{r=k}^{1+(k-1)j} |\mathcal{C}(j,k,r)|n^r p_n^{d(j,k,r)}\nonumber\\
  &\leq \frac{j^{j-1}}{a(G)^j}n^{1+(k-1)j}p_n^{je(G)}
  \sum_{r=k}^{1+(k-1)j} |\mathcal{C}(j,k,r)|\nonumber\\
  &= \frac{j^{j-1}}{a(G)^j}n^{1+(k-1)j}p_n^{je(G)}|{\rm CNF}(j,k)|\nonumber\\
  &\leq j!^kk!^{j-1} \frac{j^{j-1}}{a(G)^j}n^{1+(k-1)j}p_n^{je(G)},
\end{align*}
where the last inequality is from Lemma~\ref{fkla12}.
\item When $p_n\ll n^{-(k-1)/e(G)}$, for any
  $r\in \{k+1,\dots,(k-1)j+1 \}$, by \eqref{convex3} we have 
  \begin{align}
    \nonumber
  \frac{n^kp_n^{e(G)}}{n^rp_n^{d(j,k,r)}}&=
  \left(np_n^{\frac{d(j,k,r)-e(G)}{r-k}}\right)^{-(r-k)}
  \\
    \nonumber
  & \gtrsim\left( np_n^{\frac{e(G)}{k-1}}\right)^{-(r-k)}
  \\
  \label{fjkldfaa1-2}
  & \gg 1.
\end{align}
Hence, from \eqref{cumbound1},
 for $n$ large enough we have 
\begin{align*}
  \vert\kappa_j(N_G)\vert &\leq  \frac{j^{j-1}}{a(G)^j}
  \sum_{r=k}^{1+(k-1)j} |\mathcal{C}(j,k,r)|n^r p_n^{d(j,k,r)}\\
  &\leq \frac{j^{j-1}}{a(G)^j}n^kp_n^{e(G)}
  \sum_{r=k}^{1+(k-1)j} |\mathcal{C}(j,k,r)|\\
  &= \frac{j^{j-1}}{a(G)^j}n^kp_n^{e(G)}|{\rm CNF}(j,k)|\\
  &\leq j!^kk!^{j-1}\frac{j^{j-1}}{a(G)^j}n^kp_n^{e(G)}.
\end{align*}
\end{enumerate}
The proof is complete.
\end{Proof}
 By inspection of the proof of
  Theorem~\ref{th7.5}, see \eqref{fjkldfaa1} and \eqref{fjkldfaa1-2} therein,
  we note that its conclusions
  also hold under the following alternative conditions. 
\begin{prop}
    Let $G$ be a strongly balanced connected graph.
  \begin{enumerate}[a)]
\item If $n^{-(v(G)-1)/e(G)} \leq p_n$, $n \geq 1$, then we have 
  \begin{equation}
\nonumber %
    |\kappa_j(N_G)|\leq
   j^{j-1}j!^{v(G)} \frac{v(G)!^{j-1}}{a(G)^j}n^{1+(v(G)-1)j}p_n^{je(G)},
    \quad j \geq 2, \ n \geq 1.  
  \end{equation}
\item If $p_n \leq n^{-(v(G)-1)/e(G)}$, $n \geq 1$, then we have  
    \begin{equation}
\nonumber %
      |\kappa_j(N_G)|\leq
      j^{j-1}j!^{v(G)}
      \frac{v(G)!^{j-1}}{a(G)^j}n^{v(G)}p_n^{e(G)},
    \quad j \geq 2, \ n \geq 1.  
  \end{equation}
\end{enumerate}
\end{prop} 
 Proposition~\ref{prop7.6} is a straightforward
consequence of Proposition~\ref{p6.1}.
 Note that \eqref{var-low} in  Assumption~\ref{assu1} 
 is valid in the binomial RCM,
but it is
not satisfied in the Erd{\H o}s-R\'enyi model.    
Indeed, the principal degree of $f$ equals $2$ 
in the Erd{\H o}s-R\'enyi model,
 see \S~9 in \cite{Janson91}.
 \begin{prop} 
  \label{prop7.6}
  Let $G$ be a strongly balanced connected graph. %
  Suppose that the
  function $f$ in \eqref{fun1} satisfies Assumption~\ref{assu1}. 
\begin{enumerate}[a)]
  \item If $n^{-(v(G)-1)/e(G)}\ll p_n$, we have the lower bound 
  \begin{equation}
    \label{low-1}
    \kappa_2(N_G)\geq \frac{C}{a(G)^2}
    \frac{n!}{(n+1-2v(G))!}p_n^{2e(G)}, \quad j \geq 1, 
  \end{equation}
  where $C>0$ is a constant independent of $n \geq 2v(G)-1$.
\item If $p_n\ll n^{-(v(G)-1)/e(G)}$, we have 
  the lower bound
  \begin{equation}
    \label{low-2}
    \kappa_2(N_G)\geq \frac{C}{a(G)^2}
    \frac{n!}{(n-v(G))!}p_n^{e(G)}, \quad j \geq 1,  
  \end{equation}
  where $C>0$ is a constant  independent of $n \geq 2v(G)-1$. 
\end{enumerate}
\end{prop} 
As a consequence of Theorem~\ref{th7.5} and \ref{prop7.6}, 
we have the following corollary in the binomial RCM. 
\begin{corollary}
  \label{cjkfl} 
  Let $G$ be a strongly balanced connected graph. %
  Suppose that the function $f$ in \eqref{fun1}
  satisfies Assumption~\ref{assu1}. %
  Then, for $n\geq 1$ large enough, %
  the $j$-th cumulant of the normalized subgraph count
  $$
  \widebar{N}_G:= \frac{N_G-\kappa_1(N_G)}{\sqrt{\kappa_2(N_G)}}
  $$ 
 satisfies 
  \begin{eqnarray*}
    \kappa_j\left(\widebar{N}_G\right)\leq \begin{cases}
      \displaystyle
      \frac{j!^{1+v(G)}}{(C n)^{j/2-1}}~~~~~~~~~~{when}~~~n^{-(v(G)-1)/e(G)}\ll p_n,
      \medskip
      \\
      \displaystyle
      \frac{j!^{1+v(G)}}{\big(C n^{v(G)}p_n^{e(G)}\big)^{j/2-1}}~~~~{when}~~~p_n\ll n^{-(v(G)-1)/e(G)},
    \end{cases}
  \end{eqnarray*}
  $j\ge3$, 
  where $C$ is a constant independent of $n \geq 1$.
\end{corollary}
\begin{Proof}
  We let $k=v(G)$.
  When $n^{-(k-1)/e(G)}\ll p_n$, from the inequalities  
  \eqref{fjklf11}, \eqref{djkld1111}, 
  \eqref{upp-1} and \eqref{low-1}, we have 
  \begin{align*}
    \kappa_j\left(\widebar{N}_G\right)&= \frac{\kappa_j(N_G)}{\kappa_2(N_G)^{j/2}}\nonumber
    \\
    &\leq j!^kk!^{j-1}\frac{j^{j-1}}{a(G)^j}n^{1+(k-1)j}p_n^{je(G)}
      \left(\frac{C}{a(G)^2}
      \frac{n!}{(n-2k+1)!}p_n^{2e(G)}\right)^{-j/2}\nonumber
      \\
    &= j^{j-1}j!^kk!^{j-1}n^{1+(k-1)j}
      \left(C
      \frac{n!}{(n-2k+1)!}\right)^{-j/2}\nonumber
      \\
      &\leq j!^{k+1}e^jk!^{j-1}n^{1+(k-1)j}
        \left(C
    \left(\frac{n}{2} \right)^{2k-1}\right)^{-j/2}\nonumber\\
        &= \frac{j!^{k+1}}{n^{(j-2)/2}}e^jk!^{j-1}
          C^{-j/2}2^{j(k-1/2)}.
  \end{align*}
  Similarly, when $p_n\ll n^{-(k-1)/e(G)}$, from 
 the inequality 
$$  \frac{n!}{(n-k)!}\ge(n-k+1)^k=\left(1-\frac{k-1}n\right)^kn^k
    \ge
    \left( \frac{3n}{4} \right)^k
    $$
    and   \eqref{upp-2}-\eqref{low-2}, we have 
  \begin{align*}
    \kappa_j\left(\widebar{N}_G\right)&= \frac{\kappa_j(N_G)}{\kappa_2(N_G)^{j/2}}\nonumber\\
    &\leq j!^kk!^{j-1}\frac{j^{j-1}}{a(G)^j}n^kp_n^{e(G)}
      \left(\frac{C}{a(G)^2}
    \frac{n!}{(n-k)!}p_n^{e(G)}\right)^{-j/2}\nonumber\\
    &= j^{j-1}j!^kk!^{j-1}n^kp_n^{e(G)}
      \left(C
    \frac{n!}{(n-k)!}p_n^{e(G)}\right)^{-j/2}\nonumber\\
    &\leq \frac{(j-1)!}{\sqrt{2\pi j}}e^jj!^kk!^{j-1}n^kp_n^{e(G)}
      \left(C
    \left(\frac{3n}{4}\right)^kp_n^{e(G)}\right)^{-j/2}\nonumber\\
      &\leq e^jj!^{k+1}k!^{j-1}n^kp_n^{e(G)}
        \left(C
        \left(\frac{3n}{4} \right)^kp_n^{e(G)}\right)^{-j/2}
        \nonumber\\
        &= \frac{j!^{k+1}e^jk!^{j-1} }{\big(n^kp_n^{e(G)}\big)^{(j-2)/2}} 
                \left(C
                \left(\frac{3}{4} \right)^k \right)^{-j/2}
                .
  \end{align*}
\end{Proof}
\noindent
From Corollary~\ref{cjkfl},
we can see that the cumulants of the normalized subgraph count $\widebar{N}_G$ satisfy the Statulevi\v{c}ius growth condition 
\eqref{Statuleviciuscond2}
with $\gamma := v(G)$ and 
  \begin{eqnarray*}
 \Delta_n := \begin{cases}
      \displaystyle
      (C n)^{1/2}~~~~~~~~~~~~~~~{when}~~~n^{-(v(G)-1)/e(G)}\ll p_n,
      \medskip
      \\
      \displaystyle
      \big(C n^{v(G)}p_n^{e(G)}\big)^{1/2}~~~~{when}~~~p_n\ll n^{-(v(G)-1)/e(G)},
    \end{cases}
  \end{eqnarray*}
  where $C$ is a constant independent of $n \geq 1$.
  Therefore, from Proposition~\ref{l1}-$i)$
  we have the following results.
\begin{corollary}
\label{corg}
(Kolmogorov bound). 
  Let $G$ be a strongly balanced connected graph. %
  Suppose that
  the function $f$ in \eqref{fun1}
  satisfies Assumption~\ref{assu1}. 
 For $n\geq 1$ large enough the normalized subgraph count
  $\widebar{N}_G$ satisfies 
  \begin{eqnarray*}
    \sup_{x\in\R}\big|\IP\big(\widebar{N}_G\leq x\big)-\Phi(x)\big|\leq \begin{cases}
      \displaystyle
      \frac{C}{n^{1/(2+4v(G))}}~~~~~~~~~\quad~~~~{when}~~~n^{-(v(G)-1)/e(G)}\ll p_n,
      \medskip
      \\
      \displaystyle
      \frac{C}{\big(n^{v(G)}p_n^{e(G)}\big)^{1/(2+4v(G))}}~~~~{when}~~~p_n\ll n^{-(v(G)-1)/e(G)}. 
    \end{cases}
  \end{eqnarray*}
for some $C >0$. 
\end{corollary}
By Corollary~\ref{cjkfl} and 
Proposition~\ref{l1}-$ii)$,
 we also have the following result.
\begin{corollary}
\label{corg2}
 (Moderate deviation principle).
  Let $G$ be a strongly balanced connected graph. %
  Suppose that
  the function $f$ in \eqref{fun1}
  satisfies Assumption~\ref{assu1}. 
  Let $( a_n )_{n \geq 1}$ be a sequence of real numbers tending to infinity, and such that 
  \begin{eqnarray*}
 \begin{cases}
      \displaystyle
  a_n \ll n^{1/(2+4 v(G))} 
 ~~~~~~~~~~~~~~~~{when}~~~n^{-(v(G)-1)/e(G)}\ll p_n,
      \medskip
      \\
      \displaystyle
   a_n \ll \big(n^{v(G)}p_n^{e(G)}\big)^{1/(2+4 v(G))} 
  ~~~~{when}~~~p_n\ll n^{-(v(G)-1)/e(G)}. 
    \end{cases}
  \end{eqnarray*}
  Then, $\big(a_n^{-1}\widebar{N}_G\big)_{n \ge1}$ satisfies a moderate deviation principle with speed $a_n^2$ and rate function $x^2/2$. 
\end{corollary} 
 Another direct consequence of Theorem~\ref{th7.5} is
the following threshold phenomenon of subgraph containment in the binomial RCM.
  \begin{corollary}
    \label{graphcontain}
    Let $G$ be a strongly balanced connected graph,
    and suppose that %
    the function $f$ in \eqref{fun1} satisfies Assumption~\ref{assu1}.
    We have 
    \begin{enumerate}[a)]
      \item $\lim_{n\to\infty}\IP\left(N_G=0\right)=1$ if $p_n\ll n^{-v(G)/e(G)}$,
      \item $\lim_{n\to\infty}\IP\left(N_G=0\right)=0$ if $p_n\gg n^{-v(G)/e(G)}$.
    \end{enumerate}
  \end{corollary}
  \begin{Proof}
 We let $k=v(G)$. 
 When $p_n\gg n^{-k/e(G)}$,
 from \eqref{upp-1} we have %
 \begin{equation}
   \nonumber
   \kappa_2(N_G)\leq \frac{2^{k+1} k!}{a(G)^2} n^{1+2(k-1)}p_n^{2e(G)}, %
 \end{equation}
 hence by \eqref{mean-1} we have 
 \begin{equation}
  \nonumber
    \lim_{n\to \infty}
   \frac{\kappa_2(N_G)}{\kappa_1(N_G)^2} = 0, %
 \end{equation}
 and
 \begin{eqnarray}
 \nonumber
     \lim_{n\to \infty}
   \frac{(\E[N_G])^2}{\E[N_G^2]}
   =
   \lim_{n\to \infty} \frac{\kappa_1(N_G)^2}{\kappa_2(N_G)+\kappa_1(N_G)^2}
   = 1.
 \end{eqnarray}
 On the other hand, if $p_n\ll n^{-k/e(G)}$, from \eqref{mean-1}, we have 
 \begin{equation}
\nonumber
    \lim_{n\to \infty} \E(N_G)=0. %
 \end{equation}
 We conclude by the first and second moment methods,
 see \cite[Page~54]{JLR}, which state that 
\begin{equation}
\label{firstsecm}
\frac{(\E[N_G])^2}{\E[N_G^2]}\le\IP(N_G>0)\leq \E[N_G]. 
\end{equation}
\end{Proof}
\appendix
\section{Cumulant method} 
\label{appx}
\noindent
 Given $(X_l)_{l\geq 1}$ a sequence of random variables,
 for any subset $b$ of $\{1,\dots,n\}$ we consider a family 
\begin{equation*}
\mathbf{X}_{b}=(X_l)_{l\in b}
\end{equation*}
indexed by $b \subset [n]$.
Taking $b:=\{j_1,\dots,j_k\}$,
the joint characteristic function of the vector
$\mathbf{X}_b$ is defined as 
 $$\varphi_{\mathbf{X}_b}(t_1,\dots,t_k):=\E\left[\exp\left(i\sum_{\ell=1}^kt_\ell X_{j_\ell}\right)\right],
 $$
 and the the joint cumulant of
 $\mathbf{X}_b$
 is defined as 
 \begin{equation*}
\kappa(\mathbf{X}_b)=(-i)^k\frac{\partial^k}{\partial t_1\cdots\partial t_k}\ln \varphi_{\mathbf{X}_b}(t_1,\dots,t_k)\vert_{t_1=\cdots=t_k=0}. 
\end{equation*}
For any random variable $\xi$, we let
$\xi'_b:=(\xi,\dots,\xi)$ denote the random vector with $|b|$
entries identical to $\xi$,
and write $\kappa_n(\xi):=\kappa(\xi'_{[n]})$ for the $n$-$th$ order
 cumulant of $\xi$. 
 From Theorem~1 in \cite{elukacs},
 see also \cite{leonov}, \cite{malyshev1980}, 
 we have the relation  
 \begin{equation}
   \label{cummom1}
    \kappa(\mathbf{X}_b)=
    \sum_{\sigma
      \in\Pi(b)}(-1)^{s-1}(s-1)!\prod_{a\in \sigma}
    \E \Bigg[
      \prod_{l\in a} X_l
    \Bigg] 
  \end{equation}
 between the joint
 moments and cumulants of the random vector $\mathbf{X}_b$,
 $b\subset[n]$. 
 The following results are summarized from the ``main lemmas'' in Chapter~2 of \cite{saulis} and \cite{doring}.
\begin{prop}
\label{l1}
  Let $(X_n)_{n \geq 1}$ be a family of random variables with mean zero and unit variance for all $n\ge1$. Suppose that for all $n \geq 1$, all moments of the random variable $X_n$ exist and that %
  the cumulants of $X_n$ satisfy  
\begin{equation}
\label{Statuleviciuscond2}
    |\kappa_j (X_n)|\leq \frac{(j!)^{1+\gamma}}{(\Delta_n)^{j-2}},
    \qquad
 j\ge3, 
\end{equation}
  where $\gamma \geq 0$ depends only on $n \geq 1$,
  while $\Delta_n\in(0,\infty)$ may depend on $n \geq 1$.
   Then, the following assertions hold.
\begin{enumerate}[i)]
\item (Kolmogorov bound,
 \cite[Corollary~2.1]{saulis} and \cite[Theorem~2.4]{doering})
  One has
\begin{equation*}
\sup_{x\in\R}|\IP(X_n\leq x)-\Phi(x)|\leq \frac{C(\gamma )}{(\Delta_n)^{1/(1+2\gamma)}},
\end{equation*}
for some $C(\gamma ) >0$ depending only on $\gamma$.
\item (Moderate deviation principle,
  \cite[Theorem~1.1]{doring} and \cite[Theorem~3.1]{doering}).
  Let $( a_n )_{n \geq 1}$ be a sequence of real numbers tending to infinity, and such that 
  $$
  a_n \ll (\Delta_n)^{1/(1+2\gamma)}. 
  $$
  Then, $ (a_n^{-1}X_n)_{n \ge1}$ satisfies a moderate deviation principle with speed $a_n^2$ and rate function $x^2/2$. 
\item (Concentration inequality,
  corollary of \cite[Lemma~2.4]{saulis} and \cite[Theorem~2.5]{doering}).
  For any $x\geq 0$ and sufficiently large $n \geq 1$, 
\begin{equation*}
\IP(|X_n|\geq x)\le2\exp\left(-\frac14\min\left( \frac{x^2}{2^{1+\gamma}},(x\Delta_n)^{1/(1+\gamma)}\right) \right).
\end{equation*} 
\item (Normal approximation with Cram\'er corrections, \cite[Lemma~2.3]{saulis}). There exists a constant $c>0$ such that for all $n \geq 1$ and $x\in(0,c(\Delta_n)^{1/(1+2\gamma)})$ we have 
\begin{eqnarray*} \frac{\IP(X_n\geq x)}{1-\Phi(x)}&=& \left(1+O\left(\frac{x+1}{(\Delta_n)^{1/(1+2\gamma)}}\right)\right) \exp \big( \widetilde{L}(x) \big), \\ \frac{\IP(X_n\leq -x)}{\Phi(-x)}&=& \left(1+O\left(\frac{x+1}{(\Delta_n)^{1/(1+2\gamma)}}\right)\right) \exp \big( \widetilde{L}(-x) \big), \end{eqnarray*} where $\widetilde{L}(x)$ is related to the Cram\'er-Petrov series. 
\end{enumerate}
\end{prop}

\footnotesize

\def\cprime{$'$} \def\polhk#1{\setbox0=\hbox{#1}{\ooalign{\hidewidth
  \lower1.5ex\hbox{`}\hidewidth\crcr\unhbox0}}}
  \def\polhk#1{\setbox0=\hbox{#1}{\ooalign{\hidewidth
  \lower1.5ex\hbox{`}\hidewidth\crcr\unhbox0}}} \def\cprime{$'$}


\begin{thebibliography}{AcdOG25}

\bibitem[AcdOG25]{alvarado}
J.D. Alvarado, L.~Gon\c calves~de Oliveira, and S.~Griffiths.
\newblock Moderate deviations of triangle counts in sparse
  {E}rd{\Horig{o}}s-{R}\'enyi random graphs ${G}(n, m)$ and ${G}(n, p)$.
\newblock {\em Probab. Theory Relat. Fields}, 2025.

\bibitem[BCJ23]{bhattacharya23}
B.~Bhattacharya, A.~Chatterjee, and S.~Janson.
\newblock Fluctuations of subgraph counts in graphon based random graphs.
\newblock {\em Combin. Probab. Comput.}, 32(3):428--464, 2023.

\bibitem[BDMM24]{bhattacharya2}
B.B. Bhattacharya, S.~Das, S.~Mukherjee, and S.~Mukherjee.
\newblock Fluctuations of quadratic chaos.
\newblock {\em {Commun. Math. Phys.}}, 405(237), 2024.

\bibitem[CGR16]{coulson16}
M.~Coulson, R.~E. Gaunt, and G.~Reinert.
\newblock Poisson approximation of subgraph counts in stochastic block models
  and a graphon model.
\newblock {\em ESAIM Probab. Stat.}, 20:131--142, 2016.

\bibitem[CT22]{can2022}
V.~H. Can and K.~D. Trinh.
\newblock Random connection models in the thermodynamic regime: central limit
  theorems for add-one cost stabilizing functionals.
\newblock {\em Electron. J. Probab.}, 27:1--40, 2022.

\bibitem[DE09]{doring1}
H.~D{\"o}ring and P.~Eichelsbacher.
\newblock Moderate deviations in a random graph and for the spectrum of
  bernoulli random matrices.
\newblock {\em Electron. J. Probab.}, 14:2636--2656, 2009.

\bibitem[DE13]{doring}
H.~D{\"{o}}ring and P.~Eichelsbacher.
\newblock Moderate deviations via cumulants.
\newblock {\em J. Theoret. Probab.}, 26:360--385, 2013.

\bibitem[DF14]{DevroyeFraiman14}
L.~Devroye and N.~Fraiman.
\newblock Connectivity of inhomogeneous random graphs.
\newblock {\em Random Structures Algorithms}, 45(3):408--420, 2014.

\bibitem[DJS22]{doering}
H.~D{\"o}ring, S.~Jansen, and K.~Schubert.
\newblock The method of cumulants for the normal approximation.
\newblock {\em Probab. Surv.}, 19:185--270, 2022.

\bibitem[DM83]{mandelbaum}
E.~B. Dynkin and A.~Mandelbaum.
\newblock Symmetric statistics, {P}oisson point processes, and multiple
  {W}iener integrals.
\newblock {\em Ann. Statist.}, 11(3):739--745, 1983.

\bibitem[DST16]{dsc}
L.~Decreusefond, M.~Schulte, and C.~Th{\"a}le.
\newblock Functional {P}oisson approximation in {K}antorovich-{R}ubinstein
  distance with applications to ${U}$-statistics and stochastic geometry.
\newblock {\em Ann. Probab.}, 44(3):2147--2197, 2016.

\bibitem[ER23]{rednos}
P.~Eichelsbacher and B.~Redno{\ss}.
\newblock Kolmogorov bounds for decomposable random variables and subgraph
  counting by the {S}tein-{T}ikhomirov method.
\newblock {\em Bernoulli}, 29(3):1821--1848, 2023.

\bibitem[FMN16]{feray16}
V.~F{\'e}ray, P.-L. M{\'e}liot, and A.~Nikeghbali.
\newblock {\em Mod-$\phi$ Convergence Normality Zones and Precise Deviations}.
\newblock SpringerBriefs in Probability and Mathematical Statistics. Springer
  Cham, 2016.

\bibitem[HHO25]{heerten}
N.~Heerten, C.~Hirsch, and M.~Otto.
\newblock Moderate deviations for weight-dependent random connection models.
\newblock {\em J. Appl. Probab.}, pages 1--23, 2025.

\bibitem[Hoe61]{hoeffding61}
W.~Hoeffding.
\newblock The strong law of large numbers for {$U$}-statistics.
\newblock Technical report, North Carolina State University. Dept. of
  Statistics, 1961.

\bibitem[HP{\v{C}}21]{hladky21}
J.~Hladk{\'y}, C.~Pelekis, and M.~{\v{C}}ileikis.
\newblock A limit theorem for small cliques in inhomogeneous random graphs.
\newblock {\em J. Graph Theory}, 97:578--599, 2021.

\bibitem[Jan97]{janson}
S.~Janson.
\newblock {\em Gaussian {H}ilbert spaces}, volume 129 of {\em Cambridge Tracts
  in Mathematics}.
\newblock Cambridge University Press, Cambridge, 1997.

\bibitem[J{\L}R00]{JLR}
S.~Janson, T.~{\L}uczak, and A.~Ruci{\'n}ski.
\newblock {\em Random graphs}.
\newblock Wiley-Interscience Series in Discrete Mathematics and Optimization.
  Wiley-Interscience, New York, 2000.

\bibitem[JN91]{Janson91}
S.~Janson and K.~Nowicki.
\newblock The asymptotic distributions of generalized {$U$}-statistics with
  applications to random graphs.
\newblock {\em Probab. Theory Related Fields}, 90(3):341--375, 1991.

\bibitem[Kho08]{khorunzhiy}
O.~Khorunzhiy.
\newblock On connected diagrams and cumulants of {E}rd{\Horig{o}}s-{R}\'enyi
  matrix models.
\newblock {\em Comm. Math. Phys.}, 282:209--238, 2008.

\bibitem[KR21]{KaurRollin21}
G.~Kaur and A.~R{\"o}llin.
\newblock Higher-order fluctuations in dense random graph models.
\newblock {\em Electron. J. Probab.}, 25:1--36, 2021.

\bibitem[LNS21]{LNS21}
G.~Last, F.~Nestmann, and M.~Schulte.
\newblock The random connection model and functions of edge-marked {P}oisson
  processes: second order properties and normal approximation.
\newblock {\em Ann. Appl. Probab.}, 31(1):128--168, 2021.

\bibitem[LP24a]{LiuPrivault}
Q.~Liu and N.~Privault.
\newblock Normal approximation of subgraph counts in the random-connection
  model.
\newblock {\em Bernoulli}, 30:3224--3250, 2024.

\bibitem[LP24b]{LiuPrivault24}
Q.~Liu and N.~Privault.
\newblock Normal to {P}oisson phase transition for subgraph counting in the
  random-connection model.
\newblock Preprint arXiv:2409.16222, 37 pages, 2024.

\bibitem[{\L}R92]{luczak}
T.~{\L}uczak and A.~Ruci\'nski.
\newblock Convex hulls of dense balanced graphs.
\newblock {\em J. Comput. Appl. Math.}, 41:205--213, 1992.

\bibitem[LRP17]{lachieze-rey-peccati}
R.~Lachi\`eze-Rey and G.~Peccati.
\newblock New {B}erry-{E}sseen bounds for functionals of binomial point
  processes.
\newblock {\em Ann. Appl. Probab.}, 27(4):1992--2031, 2017.

\bibitem[LS59]{leonov}
V.P. Leonov and A.N. Shiryaev.
\newblock {On a method of calculation of semi-invariants.}
\newblock {\em Theory Probab. Appl.}, 4:319--329, 1959.

\bibitem[Luk55]{elukacs}
E.~Lukacs.
\newblock Applications of {F}a\`a di {B}runo's formula in mathematical
  statistics.
\newblock {\em Amer. Math. Monthly}, 62:340--348, 1955.

\bibitem[LX25]{liu-xia}
Q.~Liu and A.~Xia.
\newblock Poisson approximation for subgraph counts in general
  random-connection models.
\newblock Preprint, 29 pages, 2025.

\bibitem[Mal80]{malyshev1980}
V.A. Malyshev.
\newblock Cluster expansions in lattice models of statistical physics and the
  theory of quantum fields.
\newblock {\em Russian Mathematical Surveys}, 35(1):1--62, 1980.

\bibitem[MM91]{MalyshevMinlos91}
V.A. Malyshev and R.A. Minlos.
\newblock {\em Gibbs Random Fields}, volume~44 of {\em Mathematics and its
  Applications (Soviet Series)}.
\newblock Kluwer Academic Publishers Group, Dordrecht, 1991.

\bibitem[Pen03]{penrosebk}
M.~Penrose.
\newblock {\em Random Geometric Graphs}, volume~5 of {\em Oxford Studies in
  Probability}.
\newblock Oxford University Press, 2003.

\bibitem[Pen18]{penrose18}
M.~Penrose.
\newblock Inhomogeneous random graphs, isolated vertices, and {P}oisson
  approximation.
\newblock {\em J. Appl. Probab.}, 55(1):112--136, 2018.

\bibitem[PS18]{PS}
N.~Privault and G.~Serafin.
\newblock Stein approximation for functionals of independent random sequences.
\newblock {\em Electron. J. Probab.}, 23:Paper No. 4, 1--34, 2018.

\bibitem[PT11]{peccatitaqqu}
G.~Peccati and M.~Taqqu.
\newblock {\em Wiener Chaos: Moments, Cumulants and Diagrams: A survey with
  Computer Implementation}.
\newblock Bocconi \& Springer Series. Springer, 2011.

\bibitem[SS91]{saulis}
L.~Saulis and V.A. Statulevi\v{c}ius.
\newblock {\em Limit Theorems for Large Deviations}, volume~73 of {\em
  Mathematics and its Applications (Soviet Series)}.
\newblock Kluwer Academic Publishers Group, Dordrecht, 1991.

\bibitem[ST24]{schulte-thaele}
M.~Schulte and C.~Th{\"a}le.
\newblock Moderate deviations on {P}oisson chaos.
\newblock {\em Electron. J. Probab.}, 29:no. 146, 1--27, 2024.

\bibitem[Zha22]{zhangzs}
Z.S. Zhang.
\newblock Berry-{E}sseen bounds for generalized {$U$}-statistics.
\newblock {\em Electron. J. Probab.}, 27:1--36, 2022.

\end{thebibliography}
\end{document}